\documentclass[a4paper,12pt]{article}
\usepackage{amsmath}
\usepackage{color}

\newcommand{\rz}{{\rm I\!R}}
\newcommand{\nz}{{\rm I\!N}}
\newcommand{\ve}{{\varepsilon}}
\newcommand{\qed}{\hfill 
\colorbox{black}{\hspace{-0.01cm}}}

\begin{document}

\begin{center}
{\bf {\huge Distributed optimal control\\[1mm]  of a nonstandard system of \\[2mm] 
phase field equations}\\[6mm]{\small 
Dedicated to Prof. Dr. Ingo M\"uller on the occasion of\\[1mm]
  his 75th birthday}}

\vspace{9mm}
{\large Pierluigi Colli$^1\!$,
Gianni Gilardi\footnote{Dipartimento di Matematica  ``F. Casorati'',
Universit\`a di Pavia, Via Ferrata, 1, 27100 Pavia,  Italy,
e-mail: pierluigi.colli@unipv.it, gianni.gilardi@unipv.it },
Paolo Podio-Guidugli\footnote{Dipartimento di Ingegneria Civile,
Universit\`a di Roma ``Tor Vergata'',
Via del Politecnico, 1, 00133 Roma, Italy,
e-mail:ppg@uniroma2.it},
and \\[1mm]J\"urgen Sprekels\footnote{Weierstrass Institute for 
Applied Analysis and Stochastics,
Mohrenstr. 39, 10117 Berlin, Germany,
e-mail: juergen.sprekels@wias-berlin.de \\[2mm]Key words: Distributed optimal control, nonlinear phase field systems, first-order necessary optimality conditions.}%
}
\end{center}

\vspace{7mm}
\begin{abstract} \noindent
We investigate a distributed optimal control problem for a phase field 
model of Cahn-Hilliard type. The model describes two-species phase segregation  
on an atomic lattice under the presence of diffusion; it has been introduced recently in
\cite{CGPS3}, on the basis of the theory developed in \cite{PG}, and consists of a system of two
highly nonlinearly coupled PDEs. For this reason, standard arguments of optimal control theory do not apply
directly, although the control constraints and the cost functional are of standard type. 
We show that the  problem admits a solution, and we derive the first-order
necessary conditions of optimality. \end{abstract}

\section{Introduction}
Let $\Omega\subset\rz^3$ denote an open and bounded domain whose smooth
boundary $\Gamma$ has outward unit normal ${\bf n}$, let $T>0$ be a given final time, and 
let $Q:=\Omega\times (0,T)$,
$\Sigma:=\Gamma\times (0,T)$. In this paper, we study distributed 
optimal control problems of the following form:

\vspace{3mm}
{\bf (CP)} \,\,Minimize the cost functional
\begin{eqnarray}
\label{eq:1.1}
 J(u, \rho,\mu)&=&\frac 1 2 \int_\Omega |\rho(x,T)-\rho_T(x)|^2\,dx\nonumber\\
&&+ \frac{\beta_1}{2}\int_0^T\int_\Omega |\mu(x,t)-\mu_T(x,t)|^2 \,dx\,dt\nonumber\\  
&&+\,\frac{\beta_2}{2}\int_0^T\int_\Omega |u(x,t)|^2\,dx\,dt
\end{eqnarray}
subject to the state system
\begin{eqnarray}
\label{eq:1.2}
(\varepsilon+2\,\rho)\mu_t+\mu\rho_t-\Delta\mu=u\,\quad\mbox{a.\,e. in }\,Q,\\
\label{eq:1.3}
\delta\rho_t-\Delta\rho+f'(\rho)=\mu\,\quad\mbox{a.\,e. in }\,Q,\\
\label{eq:1.4}
\frac{\partial \rho}{\partial \bf n}=\frac{\partial \mu}{\partial \bf n}=0\,\quad\mbox{a.\,e. on }\Sigma,\\
\label{eq:1.5}
\rho(x,0)=\rho_0(x)\,,\quad \mu(x,0)=\mu_0(x)\,,\quad\mbox{a.\,e. in } \,\Omega,
\end{eqnarray}
and to the box control constraints
\begin{equation}
  \label{eq:1.6}
u\in U_{ad}=\left\{u\in L^\infty(Q)\,;\,0\le u\le U
\quad\,\mbox{a.\,e. in }\,Q\right\}.
\end{equation}

Here, $\beta_1\ge 0$, $\beta_2\ge 0$, $\varepsilon>0$, and $\delta>0$ are constants; 
$U\in L^\infty(Q)$ denotes a given bound, and 
$\rho_T\in L^2(\Omega)$ and $\mu_T\in L^2(Q)$ represent prescribed target functions of 
the tracking-type functional $J$. 
Although for large parts of the subsequent analysis much
more general cost functionals could be admitted, we restrict ourselves to the above situation
for the sake of a simpler exposition.

The state system (\ref{eq:1.2})--(\ref{eq:1.5}) constitutes a phase field model of Cahn-Hilliard type 
that  describes  phase segregation of two species (atoms and vacancies, say) on a lattice
in the presence of diffusion; it has been introduced recently in [15, 4]. 
The state variables  are the  {\em order parameter} $\rho$,
interpreted as a volumetric density,  and the {\em chemical potential} \,$\mu$. For physical reasons, we must 
have $0\le \rho\le 1$ and $\mu > 0$ almost everywhere in $Q$. The control
function $u$ on the right-hand side of (\ref{eq:1.2}) plays the role of a {\em microenergy source}
(see below).  Moreover, the nonlinearity $f$ is a double-well potential defined in (0,1), whose derivative 
$f'$ is singular at the endpoints $\rho=0$ and $\rho=1$: e.\,g., $f=f_1+f_2$, with
$f_2$ smooth and $f_1(\rho)=c\,(\rho\,\log(\rho)+(1-\rho)\,\log(1-\rho))$,  with $\,c\,$ a positive
constant.

System (\ref{eq:1.2})--(\ref{eq:1.5}) is singular, with highly nonlinear
and nonstandard coupling. 
In particular, nasty nonlinear terms involving time derivatives occur in (\ref{eq:1.2}),
and the expression $f'(\rho)$ in (\ref{eq:1.3}) may become singular.  
For the case $u=0$ (no control), this system was analyzed in a recent paper
\cite{CGPS3}; the case $\varepsilon\searrow 0$ was studied in \cite{CGPS4}. We also refer
to the papers \cite{CGPS1} and \cite{CGPS2}, where the corresponding Allen-Cahn model
was discussed.   

The mathematical literature on control problems for phase field systems is scarce and
usually restricted to the so-called {\em Caginalp model} of phase transitions (see, e.\,g., 
\cite{HJ}, \cite{Hei}, \cite{HeiTr}, \cite{Tr}, and the references given there). More general,
thermodynamically consistent phase field models were the subject of \cite{LS07}. Control
problems for the system (\ref{eq:1.2})--(\ref{eq:1.5}) have never been studied before.
We remark at this place that it would be a challenging task 
to study {\em boundary} control problems for the PDE system (\ref{eq:1.2}), (\ref{eq:1.3})
in place of distributed ones as in this paper; notice, however, that this would require to 
first establish appropriate well-posedness results for non-homogeneous Neumann boundary conditions or 
for non-homogeneous boundary conditions of third kind. Such results are
presently not available.

The paper is organized as follows: below, we briefly recall the thermodynamic background of the state
system (\ref{eq:1.2})--(\ref{eq:1.5}). In Section 2, we establish the existence of a solution to the
optimal control problem. First-order necessary optimality conditions, as usual given in
terms of the adjoint system and a variational inequality,  are derived in
Section 3. A large part of this analysis is devoted to proving that the control-to-state 
mapping is directionally differentiable in appropriate function spaces.  

\subsection{Some thermodynamic background}
The state equations (\ref{eq:1.2}),
(\ref{eq:1.3}) result from the
balances of microenergies and microforces postulated in a model for phase segregation and 
diffusion of atomic species on a lattice introduced in \cite{PG}, a paper we refer the reader
to for details. 
That model is a variation of the Cahn-Hilliard system
\begin{equation}
\label{eq:1.7}
 \rho_t - \kappa\, \Delta \mu =0, \,\quad  \mu= - \Delta\rho + 
f'(\rho),  
\end{equation}
when, for the sake of simplicity, the mobility coefficient $\kappa >0 $
is taken equal to one. Customarily, the equations in (\ref{eq:1.7}) are combined so as to
get the well-known \emph{Cahn-Hilliard equation}
\begin{equation}
\label{eq:1.8}
\rho_t = \kappa \,\Delta (- \Delta\rho + f'(\rho)),
\end{equation}

which describes diffusive phase separation processes in a 
two-phase material body.

A generalization of (\ref{eq:1.8}) was introduced by Fried and Gurtin in the papers
\cite{FG} and \cite{Gurtin}. Here is their line of reasoning: 

(i) \quad to regard the second equation in (\ref{eq:1.7}) 
as a {\em balance of microforces\,}:
\begin{equation}
\label{eq:1.9}
\mbox{div}\,{\boldsymbol{\xi}}+\pi+\gamma=0,
\end{equation}
where the distance microforce per unit volume 
is split into an internal part $\pi$ and an 
external part $\gamma$, and where the contact microforce per unit area of a 
surface oriented by its normal $\boldsymbol{n}$ is measured by $\boldsymbol{\xi}\cdot
\boldsymbol{n}$ 
in terms of the \emph{microstress} vector  $\boldsymbol{\xi}$;

\vspace{2mm}
(ii) \,\,\,to interpret the first equation in (\ref{eq:1.7}) as a 
{\em balance law for the order parameter\,}:
\begin{equation}
\label{eq:1.10}
\rho_t = - \mbox{div}\, \boldsymbol{h} + \sigma,
\end{equation}
where the pair $(\boldsymbol{h},  \sigma)$ is the \textit{inflow} of $\rho$; 

\vspace{2mm}
(iii) \,\,to restrict the admissible constitutive choices for $\pi,\boldsymbol{\xi}, 
\boldsymbol{h}$, and the \emph{free energy density} $\psi$, to those consistent in 
the sense of Coleman and Noll \cite{CN}, with an \emph{ad hoc} version of the Second Law 
of Thermodynamics -- namely, a postulated ``dissipation inequality that accommodates 
diffusion'' --
given in the form 
\begin{equation}
\label{eq:1.11}
\psi_t +(\pi-\mu)\,\rho_t-\boldsymbol{\xi}\cdot\nabla\rho_t+
\boldsymbol{h}\cdot\nabla\mu\leq 0
\end{equation}

(cf., in particular, Eq. (3.6) of \cite{Gurtin}). Within this framework, an admissible set of constitutive prescriptions turns out to be:
\begin{equation}
\label{eq:1.12}
\psi=\widehat\psi(\rho,\nabla\rho),\quad \widehat\pi(\rho,\nabla\rho,\mu)=\mu- \partial_\rho \widehat\psi(\rho,\nabla\rho), \quad \widehat{\boldsymbol{\xi}}(\rho,\nabla\rho)=\partial_{\nabla\rho} \widehat\psi(\rho,\nabla\rho),
\end{equation}
 together with 
\begin{equation}
\label{eq:1.13}
\boldsymbol{h} = - \boldsymbol{M}\nabla \mu , \quad \hbox{where } \, \boldsymbol{M}=\widehat{\boldsymbol{M}}(\rho,\nabla\rho,\mu, \nabla\mu).
\end{equation}
Moreover, it follows that the tensor-valued \emph{mobility mapping} $\boldsymbol{M}$ must obey the inequality
$$
\nabla \mu\cdot \widehat{\boldsymbol{M}}(\rho,\nabla\rho,\mu, \nabla\mu) \nabla\mu \geq 0 . 
$$
It follows from (\ref{eq:1.9}), (\ref{eq:1.10}), (\ref{eq:1.12}), and (\ref{eq:1.13})$_1$ that
$$\rho_t = \mbox{div}\,\left(\boldsymbol{M}\nabla\left(\partial_\rho \widehat\psi(\rho,\nabla\rho)-\mbox{div}\big(\partial_{\nabla\rho} \widehat\psi(\rho,\nabla\rho)\big)-\gamma\right)\right)+\sigma\,;
$$
the Cahn-Hilliard equation (\ref{eq:1.8}) results for the special choice
\begin{equation}
\label{eq:1.14}
\widehat\psi(\rho,\nabla\rho)= f(\rho)+\frac{1}{2}|\nabla\rho|^2,\,\quad \boldsymbol{M}
=\kappa \mathbf{1},
\end{equation}
provided that the external distance microforce $\gamma$ and the order parameter source term $\sigma$ 
are taken identically zero. 

In contrast to the theory developed by Fried and  Gurtin, the approach taken in \cite{PG} 
was the following: while step (i) was retained, the order 
parameter balance (\ref{eq:1.10}) and the dissipation inequality (\ref{eq:1.11}) were
replaced, respectively, by the \emph{microenergy balance}
\begin{equation}
\label{eq:1.15}
\varepsilon_t=e+w,\quad e:=-\mbox{div}\,{\overline {\boldsymbol{h}}}+{\overline \sigma},\quad w:=-\pi\,\rho_t+\boldsymbol{\xi}\cdot\nabla\rho_t,
\end{equation}
and the \emph{microentropy imbalance}
\begin{equation}
\label{eq:1.16}
\eta_t\geq -\mbox{div}\,\boldsymbol{h}+\sigma,\quad \boldsymbol{h}:=\mu{\overline {\boldsymbol{h}}},\quad \sigma:=\mu\,{\overline \sigma}.
\end{equation}
The salient new feature of this approach to phase segregation modeling is that the \emph{microentropy inflow} $(\boldsymbol{h},\sigma)$ is deemed proportional to the \emph{microenergy inflow} $({\overline{ \boldsymbol{h}}},{\overline\sigma})$ through the \emph{chemical potential} $\mu$; consistently, the free energy is defined to be
\begin{equation}
\label{eq:1.17}
\psi:=\varepsilon-\mu^{-1}\eta,
\end{equation}
where the chemical potential plays the same role as the \emph{coldness} in the deduction 
of the heat equation. Just as the absolute temperature is a macroscopic measure of microscopic {\em
agitation}, its inverse -- the coldness --  measures microscopic \emph{quiet}. Likewise, as argued in \cite{PG}, the chemical potential can be seen as a macroscopic 
measure of microscopic \emph{organization}; and, just as is always done for coldness, one can 
provisionally assume that $\mu$ is positive almost everywhere in $Q$. This assumption, which is
important to proving that the resulting system of field equations does have solutions, must be justified
a posteriori. The requirement that $\mu$ be positive is also the reason why we cannot admit
negative controls $u$ in the control problem (\ref{eq:1.1})--(\ref{eq:1.6}).  

Combining (\ref{eq:1.15})-(\ref{eq:1.17}), and assuming that $\mu>0$, one finds that
\begin{equation}
\label{eq:1.18}
\psi_t\leq -\eta_{}\,\partial_t (\mu^{-1})+\mu^{-1}\,{\overline {\boldsymbol{h}}}\cdot\nabla\mu-\pi\,\rho_t+\boldsymbol{\xi}\cdot\nabla\rho_t\,;
\end{equation}
this reduced dissipation inequality replaces (\ref{eq:1.11}) in filtering out
\emph{\`a la} Coleman-Noll the inadmissible constitutive choices. 

On taking all of the constitutive mappings delivering $\pi,\boldsymbol{\xi},\eta$, 
and ${\overline {\boldsymbol{h}}}$, to depend in principle on the list of variables 
$\rho,\nabla\rho,\mu,\nabla\mu$, and on choosing
\begin{equation}
\label{eq:1.19}
\psi=\widehat\psi(\rho,\nabla\rho,\mu)=-\mu\,\rho+f(\rho)+\frac{1}{2}|\nabla\rho|^2,
\end{equation}
one sees that compatibility with (\ref{eq:1.18}) implies that 
\begin{eqnarray}
\label{eq:1.20}
&& \widehat\pi(\rho,\nabla\rho,\mu)=-{\partial_{{\rho}} \widehat\psi(\rho,\nabla\rho,\mu)}=\mu-f'(\rho),\nonumber\\[1mm] 
&&\widehat{\boldsymbol{\xi}}(\rho,\nabla\rho,\mu)={\partial_{{\nabla\rho}} \widehat\psi(\rho,\nabla\rho,\mu)}=\nabla\rho, \nonumber\\[1mm] 
&&\widehat\eta(\rho,\nabla\rho,\mu)=\mu^2 \partial_{{\mu}} \widehat\psi(\rho,\nabla\rho,\mu)
=-\mu^2\rho,
\end{eqnarray}
together with
$$\widehat{\overline {\boldsymbol{h}}}(\rho,\nabla\rho,\mu,\nabla\mu) 
=  \widehat {\boldsymbol{H}}(\rho,\nabla\rho,\mu, \nabla\mu)\nabla \mu , \quad
\nabla \mu\cdot \widehat{\boldsymbol{H}}(\rho,\nabla\rho,\mu, \nabla\mu) \nabla\mu \geq 0 . 
$$
If we now choose for $\widehat{\boldsymbol{H}}$ the simplest expression 
$\boldsymbol{H}=\kappa \mathbf{1} $, implying a constant {and isotropic} mobility, and if we 
once again assume that the external distance microforce $\gamma$ and the source $\overline \sigma$ are null, then we can infer from (\ref{eq:1.20}) and (\ref{eq:1.17}) that the microforce balance 
(\ref{eq:1.9}) and the energy balance (\ref{eq:1.15}) become, respectively, 

\begin{equation}
\label{eq:1.21}
 \mbox{div}(\nabla\rho)+\mu- f'(\rho)  = 0,
\end{equation}
\begin{equation}
\label{eq:1.22}
 2\rho \,\mu_t + \mu \, \rho_t - \kappa\, \Delta\mu = 0.   
\end{equation} 

This is a nonlinear system for the unknowns $\rho$ and $\mu$, to be compared 
with system (\ref{eq:1.7}): while equations (\ref{eq:1.21}) and
(\ref{eq:1.7})$_2$ coincide, equation (\ref{eq:1.22}) is
considerably more difficult to handle than (\ref{eq:1.7})$_1$. Indeed, the latter 
 is linear while the former is not; moreover,
the time derivatives of $\rho$ and $\mu$ are both present in  (\ref{eq:1.22}),
and there are  nonconstant factors in front of both $\mu_t $ and 
$\rho_t$ that should remain positive during the entire evolution.
Note that, for nonzero microenergy source $\bar\sigma$, Eq. (\ref{eq:1.22})
becomes:
\begin{equation}
\label{eq:1.23}
2\,\rho\,\mu_t\,+\,\mu\,\rho_t-\kappa\,\Delta\mu=-\bar\sigma.
\end{equation}
In this sense, the control variable $u$ in (\ref{eq:1.2}) is nothing but $\,-\bar\sigma$. 

So far, it has not been possible to tackle the system (\ref{eq:1.21}), (\ref{eq:1.22})
(nor (\ref{eq:1.21}), (\ref{eq:1.23}))
mathematically. Not so for system (\ref{eq:1.2}), (\ref{eq:1.3}),
a regularized version of (\ref{eq:1.21}), (\ref{eq:1.23}) (with $u=-\bar\sigma$)
obtained by introducing the extra terms $\varepsilon\,\partial_t\mu $ in (\ref{eq:1.23}) and $\delta\,\partial_t\rho$ 
in (\ref{eq:1.21}), with small positive coefficients $\varepsilon$ and 
$\delta$ (our motivations for including such terms have been proposed 
and emphasized in \cite{CGPS3}). 

\section{Problem statement and existence}
\setcounter{equation}{0}
Consider the optimal control problem
(\ref{eq:1.2})--(\ref{eq:1.6}). For convenience, 
we introduce the abbreviated
notation $H=L^2(\Omega)$, $V=H^1(\Omega)$, $W=\left\{w\in H^2(\Omega)\,;\right.$
$\left.\partial w/\partial {\bf n}=0 \quad\mbox{on }\,\Gamma\right\}$. We endow these spaces with their 
standard norms, for which we use self-explaining notation like $\|\cdot\|_V$; for simplicity, we
also write $\|\cdot\|_H$ for the norm in the space $H\times H\times H$. Recall that the embeddings
$W\subset V\subset H$ are compact. Moreover, since $V$ is dense in $H$, we can identify $H$
with a subspace of $V^*$ in the usual way, i.\,e., by setting $\langle u,v\rangle_{V^*,V}
=(u,v)_H$ for all $u\in H$ and $v\in V$, 
where $\langle \cdot\,,\,\cdot\rangle_{V^*,V}$ denotes the duality pairing
between $V^*$ and $V$.
 Then also the embedding $H\subset V^*$ is compact, and
since $N\le 3$, we have the continuous Sobolev embeddings $W\subset C(\overline\Omega)$ and
$V\subset L^6(\Omega)$.

\vspace{2mm}
We make the following assumptions on the data:

\vspace{3mm}
(A1) $f=f_1+f_2$, where $f_1\in C^2(0,1)$ is convex, $f_2\in C^2[0,1]$,  and
\begin{equation}
\label{eq:2.1}
\lim_{r\searrow 0} f_1'(r)=-\infty, \quad \lim_{r\nearrow 1} f_1'(r)= +\infty.
\end{equation}
(A2) $\rho_0\in W$, $f'(\rho_0) \in H$, $\mu_0\in V\cap L^\infty(\Omega)$, and
\begin{equation}
\label{eq:2.2}
\inf\,\{\rho_0(x); \,x\in\Omega\}>0, \quad \sup\,\{\rho_0(x);\, x\in\Omega\}<1, \quad \mu_0\ge 0\,
\mbox{ a.\,e. in }\,\Omega. 
\end{equation}

\vspace{2mm}
Notice that (A2) implies that $\rho_0\in C(\overline\Omega)$, and that the convexity
of $\,f_1\,$ implies that $f(\rho_0)\in H$. 

An argumentation that parallels (and thus needs no repetition) the lines of the proofs of Theorem 2.2 and Theorem 2.3 
of \cite{CGPS3} (where we had $u=0$)  shows that the following well-posedness 
result holds for the state system (\ref{eq:1.2})--(\ref{eq:1.5}):

\vspace{5mm}
{\bf Theorem 2.1} \quad {\em Suppose that the hypotheses} (A1) {\em and} (A2) {\em are 
satisfied. Then we have:}

(i)\,\,\, {\em For every $u\in U_{ad}$ the state system} (\ref{eq:1.2})--(\ref{eq:1.5}) {\em has a unique solution}
$(\rho,\mu)$ {\em such that}
\begin{eqnarray}
\label{eq:2.3}
&&\rho\in W^{1,\infty}(0,T;H)\cap H^1(0,T;V)\cap L^\infty(0,T;W),\\[1mm]
\label{eq:2.4}
&&\mu\in H^1(0,T;H)\cap C^0([0,T];V)\cap L^2(0,T;W)\cap L^\infty(Q).
\end{eqnarray}   

(ii)\,\, {\em There are constants $0<\rho_*<\rho^*<1$, $\mu^*>0$, and $K_1^*>0$, 
depending only on the data, such that for every $u\in U_{ad}$ the 
corresponding solution $(\rho,\mu)$ satisfies}
\begin{eqnarray}
\label{eq:2.5}
& 0<\rho_*\le \rho\le\rho^*<1,\,\quad 0\le\mu\le\mu^*, 
\,\quad\mbox{a.\,e. in }\,Q,&\\[1mm]
 \label{eq:2.6}
&\|\rho\|_{W^{1,\infty}(0,T;H)\cap H^1(0,T;V)\cap L^\infty(0,T;W)}&\nonumber\\
&+\,\|\mu\|_
{H^1(0,T;H)\cap C^0([0,T];V)\cap L^2(0,T;W)\cap
  L^\infty(Q)}\le K^*_1.&
\end{eqnarray}

(iii)\, {\em Let} $\,u_1, u_2\in U_{ad}$, {\em and let} 
 $(\rho_1,\mu_1), (\rho_2,\mu_2)$ {\em be the corresponding solutions to} 
(\ref{eq:1.2})--(\ref{eq:1.5}). {\em Moreover, let}
$u=u_1-u_2$, $\rho=\rho_1-\rho_2$, $\mu=\mu_1-\mu_2$. {\em Then},
{\em for all $t\in [0,T]$,}
\begin{eqnarray}
\label{eq:2.7}
&\max\limits_{0\le s\le t}\,\left(\|\mu(s)\|_H^2\,+\,\|\rho(s)\|_V^2\right)\,
+\,{\displaystyle \int_0^t}
\left(\|\mu(s)\|_V^2\,+\,\|\rho_t(s)\|_H^2\right)\,ds\nonumber&\\[1mm]
&\le\,K_2^* {\displaystyle \int_0^t} \|u(s)\|_H^2\,ds\,,&
\end{eqnarray}
{\em with a constant$\,K_2^*>0\,$ that may depend on the data, but not on $u_1$, $u_2$.}

\vspace{5mm}
{\bf Remarks:}
1. Owing to (\ref{eq:2.7}), the solution operator $S:u\mapsto (\rho,\mu)$ is Lipschitz continuous as a
mapping from $U_{ad}$ (viewed as a subset of $L^2(Q)$) into $\left (H^1(0,T;H)\cap C^0([0,T];V)\right)\times
\left(L^2(0,T;V)\cap C^0([0,T];H)\right)$.

2. Thanks to (\ref{eq:2.5}) and to $f\in C^2(0,1)$, we have $f'(\rho)\in L^\infty(Q)$. Moreover, owing to (\ref{eq:2.4})
and to the embedding $V\subset L^6(\Omega)$, we have $\mu\in C^0([0,T];L^6(\Omega))$.  Note that (2.3)  implies, in particular, that $\rho$ is continuous from $[0,T]$ to  $H^s(\Omega)$ for all $s<2$. Now,
provided that $s$ is sufficiently large, we have $H^s(\Omega)\subset C(\overline \Omega)$; 
consequently, $\rho\in C(\overline Q)$. Hence, without loss of generality (by possibly choosing a larger $K_1^*$), we may assume that also
\begin{equation}
\label{eq:2.8}
\|\rho\|_{C(\overline Q)}\,+\,\|\mu\|_{C^0([0,T];L^6(\Omega))}\,+\,\|\rho_t\|_
{L^2(0,T;L^6(\Omega))}\, \le \,K_1^* \,.
\end{equation}

We are now prepared to prove existence for the control problem {\bf (CP)}:

\vspace{5mm}
{\bf Theorem 2.2} \quad {\em Suppose that the conditions} (A1) {\em and} (A2) {\em are satisfied. Then
the problem} {\bf (CP)} {\em has a solution} $\overline u \in U_{ad}$.

\vspace{3mm}
{\em Proof.} Let $\{u_n\}\subset U_{ad}$ be a minimizing sequence for {\bf (CP)}, and let $\{(\rho_n,
\mu_n)\}$ be the sequence of the associated solutions to (\ref{eq:1.2})--(\ref{eq:1.5}). We then can infer from (\ref{eq:2.6}) the 
existence of a triple $(\bar u,\bar\rho,\bar\mu)$ such that, for a suitable subsequence again indexed by
$n$, we have 
\begin{eqnarray*}
&&u_n\to \bar u\,\quad\mbox{weakly star in }\,L^\infty(Q),\\[1mm]
&&\rho_n\to\bar\rho\,\quad\mbox{weakly star in }\,W^{1,\infty}(0,T;H)\cap H^1(0,T;V)\cap L^\infty(0,T;W),\\[1mm]
&&\mu_n\to\bar\mu \,\quad\mbox{weakly star in }\,H^1(0,T;H)\cap L^\infty([0,T];V)\cap L^2(0,T;W).
\end{eqnarray*}
Clearly, we have that $\bar u\in U_{ad}$. Moreover, by virtue of the Aubin-Lions lemma 
(cf. [14, Thm. 5.1, p. 58]) and similar compactness results
(cf. [16, Sect. 8, Cor. 4]), we also have the strong convergences
\begin{eqnarray}
\label{eq:2.9}
&&\rho_n\to \bar\rho\,\quad\mbox{strongly in }\,C^0([0,T];H^s(\Omega)) \,\quad\mbox{for
all }\,s<2,\\[1mm]
\label{eq:2.10}
&&\mu_n\to \bar\mu\,\quad\mbox{strongly in }\,C^0([0,T];H)\cap L^2(0,T;V).
\end{eqnarray}
From this we infer, possibly selecting another subsequence again indexed by $n$, that 
$\rho_n\to\bar\rho $ pointwise a.\,e. in $Q$. In particular,
$\rho_*\le \bar\rho\le\rho^*$ a.\,e. in $Q$ and, since $f\in C^2(0,1)$, also
$f'(\rho_n)\to f'(\bar\rho)$  strongly in $L^2(Q)$. 
Now notice that the above convergences imply, in particular, that
\begin{eqnarray*}
&\hspace*{-0.6cm}&\rho_n\to \bar\rho\,\,\,\mbox{strongly in
}\,C^0([0,T];
L^6(\Omega), \,\,\partial_t\rho_n\to
\partial_t\bar\rho\,\,\,\mbox{weakly in }\,L^2(0,T;L^4(\Omega)),\\[1mm]
&\hspace*{-0.6cm}&\mu_n\to \bar\mu\,\,\,\mbox{strongly in }\,
L^2(0,T;L^4(\Omega)),\,\, \partial_t\mu_n\to
\partial_t\bar\mu\,\,\,\mbox{weakly in }\,L^2(Q).
\end{eqnarray*}
From this, it is easily verified that
\begin{eqnarray*}
&&\mu_n\,\partial_t\rho_n\to \bar\mu\,\partial_t\bar\rho\,\quad\mbox{weakly in }\,
L^1(0,T;H), \\[1mm]
&&\rho_n\,\partial_t\mu_n\to\bar\rho\,\partial_t\bar\mu\,\quad\mbox{weakly in }\,
L^2(0,T;L^{3/2}(\Omega)).
\end{eqnarray*}
In summary, if we pass to the limit as $n\to\infty$ in the state equations 
(\ref{eq:1.2})--(\ref{eq:1.5})
written for the triple $(u_n,\rho_n,\mu_n)$, we find that $(\bar\rho,\bar\mu)=S(\bar u)$,
that is, the triple $(\bar u,\bar\rho,\bar\mu)$ is admissible for the control problem
{\bf (CP)}. From the weak sequential lower semicontinuity of the cost functional $J$
it finally follows that $\bar u$, together with $(\bar\rho,\bar\mu)=S(\bar u)$, is
a solution to {\bf (CP)}. This concludes the proof. \qed

\vspace{5mm}
{\bf Remarks:} \,3. It can be shown that this existence result holds for much more general cost
functionals. All we need is that $J$ enjoy appropriate weak sequential
lower semicontinuity properties that match the above weak convergences.

\vspace{2mm}
4. Since the state component $\rho$ is continuous on $\overline Q$, the existence 
result remains valid if suitable pointwise state constraints for $\rho$
are added (provided the admissible set is not empty). For instance, consider the case
when the state has to obey the one-sided obstacle condition 
\begin{equation}
\label{eq:2.11}
\rho(x^*,t)\ge \frac 1 2\,\quad\forall\,t\in [0,T]
\end{equation}
for some fixed $x^*\in\Omega$ (which of course requires that $\rho_0(x^*)\ge 1/2$).
If the set of admissible controls in $U_{ad}$ is not empty, i.\,e., if there is at
least one $\hat u\in U_{ad}$ such that the corresponding state component $\hat\rho$ 
satisfies (\ref{eq:2.11}), then an optimal control exists. Indeed, we pick a minimizing sequence of
admissible controls $\{u_n\}\subset U_{ad}$ with associated states $(\rho_n,\mu_n)$
obeying (\ref{eq:2.11}). Since $\Omega$ is a bounded three-dimensional domain, we can infer
from (\ref{eq:2.9}) and Sobolev embeddings that $\rho_n\to\bar\rho$ uniformly in 
$\overline Q$, whence it follows that $\bar\rho$ satisfies (\ref{eq:2.11}).   

\section{Necessary optimality conditions}
\setcounter{equation}{0}

In this section, we derive the first-order necessary conditions of optimality for
problem {\bf (CP)}. In this whole section, we generally assume that the hypotheses (A1) and
(A2) are satisfied and that $\bar u\in U_{ad}$ is an optimal control
with 
associated state
$(\bar\rho,\bar\mu)$, which has the properties (\ref{eq:2.3})--(\ref{eq:2.6}) and 
(\ref{eq:2.8}); in particular, we have
$\bar\rho\in C(\overline Q)$. For technical reasons, we need to take a slightly
 smoother nonlinear term $f$; precisely, we take

\vspace{3mm}
(A3) \,\quad $f\in C^3(0,1)$.

With this assumption, we can improve the stability estimate (\ref{eq:2.7}). 
Before stating the result, let us observe that (\ref{eq:2.3}) implies, in particular,
that the solution component $\rho$ is weakly continuous from $[0,T]$ into $W$, which
justifies the formulation of the next estimate (\ref{eq:3.1}).

\vspace{7mm}
{\bf Lemma 3.1} \,\quad {\em Suppose that} (A1)--(A3) {\em are satisfied, and let
$\,u_1, u_2\in U_{ad}$ be given and 
 $(\rho_1,\mu_1), (\rho_2,\mu_2)$ be the corresponding solutions to} 
(\ref{eq:1.2})--(\ref{eq:1.5}). {\em Moreover, let}
$u=u_1-u_2$, $\rho=\rho_1-\rho_2$, $\mu=\mu_1-\mu_2$. {\em Then},
{\em for all $t\in [0,T]$,}
\begin{eqnarray}
\label{eq:3.1}
\max_{0\le s\le t}\,\left(\|\rho_t(s)\|_V^2 \,+\,\|\mu(s)\|_V^2\,
+\,\|\rho(s)\|_W^2\right)
\,+\,\int_0^t\!\!
\left(\|\mu_t(s)\|_H^2\,+\,\|\rho_t(s)\|_W^2\right)\,ds\nonumber\\
\,\le\,K_3^*
\int_0^t\|u(s)\|_H^2\,ds\,,\hspace{3cm}
\end{eqnarray}
{\em with a constant$\,K_3^*>0\,$ that may depend on the data, but not on $u_1$, $u_2$.}

\vspace{2mm}
{\em Proof.} \,\quad Obviously, the pair $(\rho,\mu)$ is a solution to the system 
\begin{eqnarray}
\label{eq:3.2}
(\ve+2\rho_1)\mu_t \,+\,2\,\rho\,\mu_{2,t}\,+\,\mu\,\rho_{1,t}+\mu_2\,\rho_t-\Delta
\mu=u\,\quad\mbox{a.\,e. in }\,Q,\\[1mm]
\label{eq:3.3}
\delta\,\rho_t-\Delta\rho=\mu\,-\,(f'(\rho_1)-f'(\rho_2))\,\quad\mbox{a.\,e. in }\,Q,\\[1mm]
\label{eq:3.4}
\frac{\partial \rho}{\partial \bf n}=\frac{\partial \mu}{\partial \bf n}=0\,\quad\mbox{a.\,e. on }\Sigma,\\[1mm]
\label{eq:3.5}
\rho(x,0)= \mu(x,0)=0\,,\quad\mbox{a.\,e. in } \,\Omega.
\end{eqnarray}
We test Eq. (\ref{eq:3.2}) by $\mu_t$. It then follows, with the use of Young's inequality, that
\begin{eqnarray}
\label{eq:3.6}
\frac{\ve} 2\int_0^t\|\mu_t(s)\|^2_H\,ds \,+\,\frac 1 2 \|\nabla\mu(t)\|_H^2\,\le\,\frac 1 2
\int_0^t\|u(s)\|_H^2\,ds\nonumber\\
+ \int_0^t\int_\Omega\left(2|\rho|\,|\mu_{2,t}|\,+\,|\mu|\,|\rho_{1,t}|\,+\,|\mu_2|
\,|\rho_t|\right)|\mu_t|\,dx\,ds\,.
\end{eqnarray}

We estimate the terms on the right-hand side individually. In this process, $C_i$
($i\in\nz$)  denote positive constants that only depend on the constants $\,\varepsilon, \delta,
\rho_*,\rho^*,\mu^*,T,K_1^*,K_2^*$.
On using H\"older's and Young's inequalities, as well as the continuity of the
embedding $H^1(\Omega)\subset L^4(\Omega)$, we have that
\begin{eqnarray}
\label{eq:3.7}
\int_0^t\!\!\int_\Omega|\mu|\,|\rho_{1,t}|\,|\mu_t|\,dx\,ds \le\int_0^t\|\mu_t(s)\|_H\,\|\mu(s)\|_{L^4(\Omega)}
\,\|\rho_{1,t}(s)\|_{L^4(\Omega)}\,ds\nonumber\\
\le \gamma \int_0^t\|\mu_t(s)\|_H^2\,ds
\,+\,\frac{C_1}{\gamma}\,\int_0^t\|\rho_{1,t}(s)\|_V^2\|\mu(s)\|_V^2\,ds\,.
\end{eqnarray}

Observe that, owing to (\ref{eq:2.6}), the function $\,s\mapsto\|\rho_{1,t}(s)\|_V^2$ belongs to $L^1(0,T)$.
Moreover, by (\ref{eq:2.5}) and Young's inequality,
\begin{equation}
\int_0^t\!\!\int_\Omega |\mu_2|\,|\rho_t|\,|\mu_t|\,dx\,ds
\,\le\,\gamma\int_0^t\|\mu_t(s)\|_H^2\,ds\,+\,\frac{C_2}{\gamma}\int_0^t\|\rho_t(s)|_H^2\,ds\,,
\end{equation}

where the second integral on the 
right-hand side can be estimated using (\ref{eq:2.7}).
In addition, we have
\begin{eqnarray}
\label{eq:3.9}
&{\displaystyle I_1:=\int_0^t\!\!\int_\Omega
\!\!2\,|\rho|\,|\mu_{2,t}|\,|\mu_t|\,dx\,ds}&\nonumber\\[1mm]
&{\displaystyle \le\gamma\int_0^t\!\!
\|\mu_t(s)\|_H^2\,ds+
\frac{C_3}{\gamma}\int_0^t\!\!\|\mu_{2,t}(s)\|_H^2\,\|\rho(s)
\|^2_{L^\infty(\Omega)}\,ds\,.}&
\end{eqnarray}
Now, observe that the embedding $W\subset L^\infty(\Omega)$, in combination with
standard elliptic estimates and (\ref{eq:2.7}),
 implies that
\begin{equation}
\label{eq:3.10}
\|\rho(s)\|^2_{L^\infty(\Omega)}\,\le\,C_4\int_0^s\|u(\tau)\|_H^2\,d\tau\,+\,C_5\,\|\Delta\rho(s)\|^2_H,
\end{equation}
whence 
\begin{eqnarray}
\label{eq:3.11}
&{\displaystyle I_1\,\le\,\gamma\int_0^t\!\!\|\mu_t(s)\|_H^2\,ds\,}&\nonumber\\[1mm]
&{\displaystyle +\,\frac{C_6}{\gamma}\Bigl(\int_0^t\!\!\|u(s)\|_H^2\,ds\,+\,
\int_0^t\!\!\|\mu_{2,t}(s)\|_H^2\,\|\Delta\rho(s)\|_H^2\,ds\Bigr)\,,}&
\end{eqnarray}
where the mapping $\,s\mapsto \|\mu_{2,t}(s)\|_H^2\,$ belongs to $L^1(0,T)$.

Next, we formally test Eq. (\ref{eq:3.3}) by $\,-\Delta\rho_t$. On integrating by parts, we find that
\begin{equation}
\label{eq:3.12}
\delta\int_0^t \!\!\|\nabla\rho_t(s)\|_H^2\,ds\,+\,\frac 1 2\, \|\Delta\rho(t)\|^2_H
\,\le\,\int_0^t\!\!\int_\Omega\!\!\left(-\mu+(f'(\rho_1)-f'(\rho_2)\right)\,\Delta\rho_t\,dx\,ds\,.\quad
\end{equation}
After a further integration by parts, this time with respect to $t$, and invoking Young's inequality and (\ref{eq:2.7}), 
we find that
\begin{eqnarray}
\label{eq:3.13}
&\hspace*{-5mm}&-\int_0^t\!\!\int_\Omega\mu\,\Delta\rho_t\,dx\,ds
=-\int_\Omega \mu(t)\,\Delta\rho(t)\,dx+
\int_0^t\!\!\int_\Omega \mu_t\,\Delta\rho\,dx\,ds\le\frac 1 8\,
\|\Delta\rho(t)\|_H^2\,
\nonumber\\[2mm]
&\hspace*{-5mm}&\quad+\,\gamma\int_0^t\!\!\|\mu_t(s)\|_H^2\,ds
+\frac{C_7}{\gamma}\!\int_0^t\!\!\|\Delta\rho(s)\|_H^2\,ds\,+\,C_8\!\int_0^t\!\!\|u(s)\|_H^2\,ds\,.
\end{eqnarray}

\noindent
Moreover, integration by parts with respect to $t$, with the help of Young's inequality, (\ref{eq:2.6}), 
and (\ref{eq:2.7}),  yields:
\begin{eqnarray}
\label{eq:3.14}
&&\int_0^t\!\!\int_\Omega (f'(\rho_1)-f'(\rho_2))\,\Delta\rho_t\,dx\,ds\,\le\,\int_\Omega |f'(\rho_1(t))
-f'(\rho_2(t))|\,|\Delta\rho(t)|\,dx\nonumber\\[1mm]
&&\quad +\,\int_0^t\!\!\int_\Omega\!\!
\left(|f''(\rho_1)|\,|\rho_t|\,+\,|f''(\rho_1)-f''(\rho_2)|\,|\rho_{2,t}|\right)\,
|\Delta\rho|\,dx\,ds\nonumber\\[2mm]
&&\le \,\frac 1 8\,\|\Delta\rho(t)\|_H^2 \,+\,\int_0^t\!\!\|\Delta\rho(s)\|_H^2\,ds \,+\,
C_9 \int_0^t\!\!\|u(s)\|_H^2\,ds\,+\,I_2\,,
\end{eqnarray}
where
$$I_2:=\int_0^t\!\!\int_\Omega\!\!|f''(\rho_1)-f''(\rho_2)|\,|\rho_{2,t}|\,
|\Delta\rho|\,dx\,ds\,.
$$

From this inequality, on applying the mean value theorem, (\ref{eq:2.6}), (\ref{eq:2.7}), Young's inequality, 
and the continuity of the
embedding $H^1(\Omega)\subset L^4(\Omega)$, we deduce the following estimate:
\begin{eqnarray}
\label{eq:3.15}
I_2&\!\!\le\!\!&C_{10}\int_0^t\!\!\|\Delta\rho(s)\|_H\,\|\rho(s)\|_{L^4(\Omega)}\,\|\rho_{2,t}(s)\|_{L^4(\Omega)}\,ds
\nonumber\\[1mm]
&\!\!\le\!\!& C_{11}\,\max_{0\le s\le t}\,\|\rho(s)\|_V\,\int_0^t\!\|\rho_{2,t}(s)\|_V\,\|\Delta\rho(s)\|_H\,ds
\nonumber\\[1mm]
&\!\!\le\!\!&C_{12}\,\Bigl(\int_0^t\!\|\rho_{2,t}(s)\|_V^2\,\|\Delta\rho(s)\|_H^2\,ds\,+\,
\int_0^t\!\|u(s)\|_H^2\,ds\Bigr)\,.
\end{eqnarray} 
Observe that by (\ref{eq:2.6}) the mapping $\,s\mapsto \|\rho_{2,t}(s)\|_V^2\,$ belongs to $L^1(0,T)$. 

At this point, we may combine the estimates (\ref{eq:3.6})--(\ref{eq:3.15}): in fact, choosing $\gamma>0$ appropriately small,
and invoking Gronwall's lemma, we find the estimate:
\begin{eqnarray}
\label{eq:3.16}
&&\int_0^t\!\! \Bigl(\|\mu_t(s)\|^2_H\,+\,\|\nabla\rho_t(s)\|_H^2\Bigr)\,ds
\,+\,\max_{0\le s \le t}\,\left(\|\mu(s)\|^2_V\,+\,\|\rho(s)\|_W^2\right)\nonumber\\[1mm]
&&\quad\le\, C_{13}\int_0^t\!\|u(s)\|_H^2\,ds\,,
\end{eqnarray}

for all $t\in [0,T]$.
Next, we formally differentiate Eq. (\ref{eq:3.3}) with respect to $t$, and obtain
\begin{equation}
\label{eq:3.17}
\delta\rho_{tt}-\Delta\rho_t=\mu_t\,-\,f''(\rho_1)\,\rho_t -(f''(\rho_1)-f''(\rho_2))\,\rho_{2,t}\,,
\end{equation}
with zero initial and Neumann boundary conditions for $\rho_t$ (cf. (\ref{eq:3.4}), (\ref{eq:3.5}),
and (\ref{eq:1.5})). Hence, testing 
(\ref{eq:3.17}) by $\rho_t$, invoking
Young's inequality, and recalling (\ref{eq:2.7}), and (\ref{eq:3.16}), we find that 
\begin{eqnarray}
\label{eq:3.18}
\frac {\delta} {2}\, \|\rho_t(t)\|_H^2\,+\,\int_0^t\!\!\|\nabla\rho_t(s)\|_H^2\,ds
\,\le\,C_{14}\int_0^t\!\!\|u(s)\|_H^2\,ds
\nonumber\\
+\,\int_0^t\!\!\int_\Omega|\rho_{2,t}|\,|f''(\rho_1)-f''(\rho_2)|\,|\rho_t|\,dx\,ds\,. 
\end{eqnarray}

Moreover, using H\"older's and Young's inequalities, (A3), (\ref{eq:2.6}), and (\ref{eq:2.7}), 
we see that

\begin{eqnarray}
\label{eq:3.19}
&&\int_0^t\!\!\int_\Omega|\rho_{2,t}|\,|f''(\rho_1)-f''(\rho_2)|\,|\rho_t|\,dx\,ds\nonumber\\
&&\,\le\,C_{15}\int_0^t\|\rho_{2,t}(s)\|_{L^4(\Omega)}\,\|\rho(s)\|_{L^4(\Omega)}\|\rho_t(s)\|_H\,ds\nonumber\\
&&\,\le\,C_{16}\,\Bigl(\int_0^t\|\rho_t(s)\|_H^2\,ds\,+\,\max_{0\le s\le t}\|\rho(s)\|_V^2
\int_0^t\|\rho_{2,t}(s)\|_V^2\,ds\Bigr)\nonumber\\
&&\,\le\,C_{17}\int_0^t\|u(s)\|_H^2\,ds\,.
\end{eqnarray}

\vspace{2mm}
Finally, we test (\ref{eq:3.17}) by $\,-\Delta\rho_t$. Using Young's inequality and
(\ref{eq:3.16}), we find that

\begin{eqnarray}
\label{eq:3.20}
&&\frac {\delta} 2 \,\|\nabla\rho_t(t)\|_H^2\,+\, \int_0^t\|\Delta\rho_t(s)\|_H^2\,ds
\,\le\,\gamma\int_0^t\|\Delta\rho_t(s)\|_H^2\,ds\,\nonumber\\[1mm]
&&\quad +\,\frac{C_{18}}{\gamma}\int_0^t\|u(s)\|_H^2\,ds\,
+\,\int_0^t\!\!\int_\Omega|\rho_{2,t}|\,|f''(\rho_1)-f''(\rho_2)|\,|\Delta\rho_t|\,dx\,ds
\nonumber\\[1mm]
&&\le \,2\gamma\!\int_0^t\!\!\|\Delta\rho_t(s)\|_H^2\,ds\,\nonumber\\[1mm]
&&\quad +\,\frac{C_{19}}{\gamma}
\Bigl(\int_0^t\!\!\|u(s)\|_H^2\,ds\,+\,\,\max_{0\le s\le t}\|\rho(s)\|_V^2
\!\int_0^t\!\!\|\rho_{2,t}(s)\|_V^2\,ds\Bigr)\nonumber\\[1mm]
&&\le \,2\gamma\!\int_0^t\!\!\|\Delta\rho_t(s)\|_H^2\,ds\,+\,\frac{C_{20}}{\gamma}
\int_0^t\!\!\|u(s)\|_H^2\,ds\,.
\end{eqnarray}

Choosing $\gamma>0$ appropriately small, we can infer that the estimate (\ref{eq:3.1}) is 
in fact true. This concludes the proof.
\qed

\subsection{The linearized system}

Suppose that $h\in L^\infty(Q)$ is an admissible variation with respect to $\bar u$, i.\,e., that there exists $\bar\lambda>0$ such that $\bar u+\lambda h\in U_{ad}$ whenever $0<\lambda\le\bar\lambda$. We have
to determine the directional derivative $DJ_{\rm red}(\bar u)h$ of the ``reduced'' cost functional
$J_{\rm red}(u):=J(u, Su)$ at $\bar u$ in the direction $h$. This requires to find the
directional derivative $DS(\bar u)h$ of the solution operator $S$ at $\bar u$ in the direction $h$. To this end, we consider the following system, which is obtained by linearizing the system 
(\ref{eq:1.2})--(\ref{eq:1.5}) at $(\bar\rho,\bar\mu)$:
\begin{eqnarray}
\label{eq:3.21}
(\ve+2\bar\rho)\,\eta_t-\Delta\eta+2\,\bar\mu_t\,\xi+\bar\mu\,\xi_t+\bar\rho_t\,\eta=h\,\quad
\mbox{a.\,e. in }\,Q,\\[1mm]
\label{eq:3.22}
\delta\,\xi_t-\Delta\xi=-f''(\bar\rho)\,\xi\,+\,\eta\,\quad
\mbox{a.\,e. in }\,Q,\\[1mm]
\label{eq:3.23}
\,\,\frac{\partial\xi}{\partial \bf n}=\frac{\partial\eta}{\partial \bf n}=0\,\quad
\mbox{a.\,e. on }\Sigma,\\[1mm]
\label{eq:3.24}
\xi(x,0)=\eta(x,0)=0 \,\,\quad \mbox{a.\,e. in }\,\Omega.
\end{eqnarray}

We expect that $(\xi,\eta)=DS(\bar u)h$, provided that (\ref{eq:3.21})--(\ref{eq:3.24}) admits a unique solution
$(\xi,\eta)$. In view of (\ref{eq:2.3}) and (\ref{eq:2.4}), we can guess the regularity of $\,\xi\,$ and $\,\eta\,$:
\begin{eqnarray}
\label{eq:3.25}
\xi\in H^1(0,T;H)\cap C^0([0,T];V)\cap L^2(0,T;W)\cap L^\infty(Q),\\
\label{eq:3.26}
\eta\in  H^1(0,T;H)\cap C^0([0,T];V)\cap L^2(0,T;W).
\end{eqnarray}

Indeed, if (\ref{eq:3.25}) and (\ref{eq:3.26}) hold then the collection of source terms in 
(\ref{eq:3.21}), i.\,e., the part $\,h-2\,\bar\mu_t\,\xi-\bar\mu\,\xi_t-\bar\rho_t\,\eta$,
belongs to $L^2(Q)$ (as it should for a solution $\eta$ satisfying (\ref{eq:3.26})),
whereas the regularity (\ref{eq:3.26}) for $\eta$ allows us to conclude from (\ref{eq:3.22})
that also $\xi\in C(\overline Q)$ (by applying maximal parabolic regularity theory, see, e.\,g.,
[7, Thm. 6.8] or [17, Lemma 7.12]).

In fact, as to $\xi$, we can count on an even better regularity. Indeed, we may differentiate (\ref{eq:3.22}) with respect to $t$ to find that
\begin{equation} 
\label{eq:3.27}
\delta\xi_{tt}-\Delta\xi_t=-f'''(\bar\rho)\,\bar\rho_t\,\xi-f''(\bar\rho)\,\xi_t+\eta_t\,,
\end{equation}
with zero initial and Neumann boundary conditions for $\,\xi_t$. Since the right-hand side of (\ref{eq:3.27})
belongs to $L^2(Q)$, we may test by any of the functions
$\xi_t$, $\xi_{tt}$, and $\,-\Delta\xi_t$, to obtain that even
\begin{equation}
\label{eq:3.28}
\xi\in H^2(0,T;H)\cap C^1([0,T];V)\cap H^1(0,T;W)\,.
\end{equation} 
Notice, however, that this fact
 has no bearing on the regularity of $\eta$, since the coefficient $\bar\mu_t$ in (\ref{eq:3.21})
only belongs to $L^2(Q)$.  

\vspace{3mm} We first prove the well-posedness of the linear system (\ref{eq:3.21})--(\ref{eq:3.24}).

\vspace{3mm} 
{\bf Proposition 3.2} \quad {\em Suppose that} (A1)--(A3) {\em are fulfilled. Then
 the system} (\ref{eq:3.21})--(\ref{eq:3.24})
 {\em has a unique solution 
$(\xi,\eta)$ satisfying} (\ref{eq:3.26}) {\em and} (\ref{eq:3.28}).     
     
\vspace{2mm}
{\em Proof.} \,We proceed in series of steps. 

\underline{Step 1: \,Approximation.} \,\quad Following the lines of our approach in \cite{CGPS3}, we use an approximation technique based on a delay in the right-hand side of (\ref{eq:3.22}). To this
end, we define for $\tau>0$ the translation operator ${\cal T}_\tau:\,L^1(0,T;H)$ $\to L^1(0,T;H)$ by putting, for every $v\in L^1(0,T;H)$ and almost every $t\in (0,T)$, 
\begin{equation}
\label{eq:3.29}
({\cal T}_\tau v)(t)= v(t-\tau)\,\mbox{ if }\,t\ge\tau, \quad\,\mbox{and }\,\,\,
({\cal T}_\tau v)(t)= 0\,\mbox{ if }\,t<\tau. 
\end{equation}
Notice that, for any $v\in L^2(Q)$ and any $\tau>0$, we obviously have 
$\,\|{\cal T}_\tau v\|_{L^2(Q)}$ $\le \|v\|_{L^2(Q)}$.

Then, for any fixed $\tau>0$, we look for functions $(\xi^\tau, \eta^\tau)$, which satisfy (\ref{eq:3.25}) and (\ref{eq:3.26}) and the system: 
\begin{eqnarray}
\label{eq:3.30}
(\ve+2\bar\rho)\,\eta^\tau_t-\Delta\eta^\tau+2\,\bar\mu_t\,\xi^\tau
+\bar\mu\,\xi^\tau_t+\bar\rho_t\,\eta^\tau=h\,\quad
\mbox{a.\,e. in }\,Q,\\[1mm]
\label{eq:3.31}
\delta\,\xi^\tau_t-\Delta\xi^\tau+f''(\bar\rho)\,\xi^\tau={\cal T}_\tau\eta^\tau\,\quad
\mbox{a.\,e. in }\,Q,\\[1mm]
\label{eq:3.32}
\,\,\frac{\partial\xi^\tau}{\partial \bf n}=\frac{\partial\eta^\tau}{\partial \bf n}=0\,\quad
\mbox{a.\,e. on }\Sigma,\\[1mm]
\label{eq:3.33}
\xi^\tau(x,0)=\eta^\tau(x,0)=0 \,\,\quad \mbox{a.\,e. in }\,\Omega.
\end{eqnarray}
Precisely, we choose for $\tau>0$ the discrete values $\tau=T/N$, where $N\in\nz$ is arbitrary, and put
$t_n=n\,\tau$, $0\le n\le N$, and  $I_n=(0,t_n)$. For $1\le n\le N$, we solve the problem
\begin{eqnarray}
\label{eq:3.34}
(\ve+2\bar\rho)\,\eta_{n,t}-\Delta\eta_n+2\,\bar\mu_t\,\xi_n+\bar\mu\,\xi_{n,t}+\bar\rho_t\,\eta_n=h\,\quad
\mbox{a.\,e. in }\,\Omega\times I_n,\\[1mm]
\label{eq:3.35}
\frac{\partial\eta_n}{\partial \bf n}=0\,\quad\mbox{a.\,e. on }\,\Gamma\times I_n, \,\quad\quad
\eta_n(x,0)=0 \,\,\quad \mbox{a.\,e. in }\,\Omega, \\[1mm]
\label{eq:3.36}
\delta\,\xi_{n,t}-\Delta\xi_n+f''(\bar\rho)\,\xi_n={\cal T}_\tau\eta_n\,\quad
\mbox{a.\,e. in }\,\Omega\times I_n,\\[1mm]
\label{eq:3.37}
\frac{\partial\xi_n}{\partial \bf n}= 0 \,\quad\mbox{a.\,e. on }\,\Gamma\times I_n,\,\quad\quad
\xi_n(x,0)=0 \,\,\quad \mbox{a.\,e. in }\,\Omega,
\end{eqnarray}

where the variables $\eta_n$ and $\xi_n$, defined on $I_n$, have obvious meaning.
Here, ${\cal T}_\tau$ acts on functions that are not defined on the entire interval $(0,T)$; however,
for $n>1$ it is still defined by (\ref{eq:3.29}), while for $n=1$ we simply 
put ${\cal T}_\tau\eta_n=0$.
Notice that whenever the pairs $(\xi_k,\eta_k)$ with
\begin{eqnarray}
\label{eq:3.38}
\xi_k\in H^1(I_k;H)\cap C^0(\bar I_k;V)\cap L^2(I_k;W)\cap C(\overline {\Omega\times I_k}),\\
\label{eq:3.39}
\eta_k\in  H^1(I_k;H)\cap C^0(\bar I_k;V)\cap L^2(I_k;W),
\end{eqnarray}
have been constructed for $1\le k\le n<N$, then we look for the pair $(\xi_{n+1},\eta_{n+1})$
that coincides with $(\xi_n,\eta_n)$ in $I_n$, and note that the linear parabolic problem (\ref{eq:3.36}), (\ref{eq:3.37}) has a unique solution $\xi_{n+1}$ on $\Omega\times I_{n+1}$ that satisfies (\ref{eq:3.38}) for $k=n+1$. Inserting
$\xi_{n+1}$ in (\ref{eq:3.34}) 
(where $n$ is replaced by $n+1$), we then find that the linear parabolic problem 
(\ref{eq:3.34}), (\ref{eq:3.35}) admits a
unique solution $\eta_{n+1}$ that fulfills (\ref{eq:3.39}) for $k=n+1$. Hence, we conclude that 
$(\xi^\tau ,\eta^\tau)=(\xi_N,\eta_N)$ satisfies (\ref{eq:3.30})--(\ref{eq:3.33}), and (\ref{eq:3.25}), 
(\ref{eq:3.26}).  

\vspace{3mm}
\underline{Step 2: \,\,A priori estimates.} \,\quad We now prove a series of a priori estimates for the functions $(\xi^\tau,\eta^\tau)$. In the following, we denote by $C_i$ ($i\in\nz$) some generic positive constants, which may depend on $\,\varepsilon, \delta,\rho_*, \rho^*,\mu^*, T, K^*_1, K^*_2$, but not on $\tau$ (i.\,e., not on $N$). For the sake of simplicity,
we omit the superscript $\tau$ and simply write $(\xi,\eta)$. We recall the continuity of the embedding
$H^1(\Omega)\subset L^6(\Omega)$. 

\vspace{2mm}
\underline{First a priori estimate.} \,\quad Observe that $\,2\,\bar\rho\,\eta\,\eta_t=
\left(\bar\rho\,\eta^2\right)_t-\bar\rho_t\,\eta^2$. Hence, testing (\ref{eq:3.30}) by $\eta$, we have, for 
$0\le t\le T$,  
\begin{eqnarray}
\label{eq:3.40}
&&\int_\Omega \left(\frac {\ve}{2}+\bar\rho\right)\eta(t)^2\,dx + \int_0^t \|\nabla\eta(s)\|_H^2\,ds
\le \frac 1 2\int_0^t \left(\|\eta(s)\|_H^2\,+\,\|h(s)\|_H^2\right)\,ds\nonumber\\
&&\quad+\,2\int_0^t\int_\Omega
|\bar\mu_t|\,|\xi|\,|\eta|\,dx\,ds \,+\int_0^t\int_\Omega
|\bar\mu|\,|\xi_t|\,|\eta|\,dx\,ds \,.
\end{eqnarray}
For any $\gamma>0$, we have by Young's inequality that
\begin{eqnarray}
\label{eq:3.41}
&&\int_0^t\int_\Omega
|\bar\mu|\,|\xi_t|\,|\eta|\,dx\,ds \,\le\,\|\bar\mu\|_{L^\infty(Q)}
\int_0^t\|\eta(s)\|_H\,\|\xi_t(s)\|_H\,ds \nonumber\\
&&\quad\le\,\gamma \int_0^t\|\xi_t(s)\|_H^2\,ds\,+\,\frac {C_1}{\gamma}\int_0^t\|\eta(s)\|_H^2\,ds\,.
\end{eqnarray}
Moreover,
\begin{eqnarray}
\label{eq:3.42}
&&\int_0^t\int_\Omega |\bar\mu_t|\,|\xi|\,|\eta|\,dx\,ds\le \int_0^t\|\bar\mu_t(s)\|_H\,\|\xi(s)\|_
{L^4(\Omega)}\,\|\eta(s)\|_{L^4(\Omega)}\,ds\nonumber\\
&&\quad \le\,\gamma\int_0^t\|\eta(s)\|^2_V\,ds\,+\,\frac{C_2}{\gamma}\,\int_0^t\|\bar\mu_t(s)\|_H^2\,
\|\xi(s)\|^2_V\,ds\,.
\end{eqnarray}
Notice that, by virtue of (\ref{eq:2.6}), the mapping $\,s\mapsto \|\bar\mu_t(s)\|_H^2\,$ belongs to
$L^1(0,T)$.

Next, we add $\xi$ on both sides of Eq. (\ref{eq:3.31}) and test the resulting equation by $\xi_t$. On
using Young's inequality again, we obtain:
\begin{eqnarray}
\label{eq:3.43}
&&\frac {\delta} 4 \int_0^t\|\xi_t(s)\|_H^2\,ds \,+\,\frac 1 2 \left(\|\xi(t)\|_H^2\,+\,\|\nabla\xi(t)\|_H^2
\right)\nonumber\\
&&\quad\le \,C_3\,\Bigl(\int_0^t\|\eta(s)\|_H^2\,ds \,+\,\int_0^t\|\xi(s)\|_H^2\,ds\Bigr)\,.
\end{eqnarray} 
Adding the inequalities (\ref{eq:3.40}) and (\ref{eq:3.43}), and
choosing $\gamma>0$ sufficiently small, we conclude from the above
estimates and Gronwall's lemma that 
\begin{eqnarray}
\label{eq:3.44}
&{\displaystyle \int_0^T\left(\|\xi_t(t)\|_H^2\,+\,\|\eta(t)\|_V^2\right)dt\,
+\,\max_{0\le t\le T}\left(
\|\xi(t)\|_V^2\,+\,\|\eta(t)\|_H^2\right)}&\nonumber\\[1mm]
&{\displaystyle \le
C_4\int_0^T\|h(t)\|_H^2\,dt.}&
\end{eqnarray}
By comparison in (\ref{eq:3.31}), and thanks to (\ref{eq:3.32}),
we may also infer (possibly by choosing a larger $C_4$) that
\begin{equation}
\label{eq:3.45}
\int_0^T\|\xi(t)\|_W^2\,dt\,\le\,C_4\int_0^T\|h(t)\|_H^2\,dt\,.
\end{equation}

\vspace{5mm}
\noindent
\underline{Second a priori estimate.} \,\quad We test (\ref{eq:3.30}) by $\eta_t$ and apply Young's
inequality in order to obtain 
\begin{eqnarray}
\label{eq:3.46}
&\hspace*{-15mm}&\frac {\ve} {2}\int_0^t\|\eta_t(s)\|_H^2\,ds 
\,+\,\frac  1 2\,\|\nabla\eta(t)\|_H^2\,\nonumber\\[1mm]
&\hspace*{-15mm}&\quad\le\,
\frac 1 {2\ve}\int_0^t\|h(s)\|_H^2\,ds+\int_0^t\!\!\int_\Omega\left(2\,|\bar\mu_t|\,|\xi|
+|\bar\mu|\,|\xi_t|+|\bar\rho_t|\,|\eta|\,\right)\,|\eta_t|\,dx\,ds.\quad
\end{eqnarray}
Since $\bar\mu\in L^\infty(Q)$, we can infer from Young's inequality that
\begin{equation}
\label{eq:3.47}
\int_0^t\int_\Omega |\bar\mu|\,|\xi_t|\,|\eta_t|\,dx\,ds \,\le\,
\gamma\int_0^t\|\eta_t(s)\|_H^2\,ds\,+\,\frac{C_5}{\gamma}\int_0^t\|\xi_t(s)\|_H^2\,ds\,.
\end{equation}
Moreover, by virtue of H\"older's and Young's inequalities,
\begin{eqnarray}
\label{eq:3.48}
&\hspace*{-6mm}&\int_0^t\int_\Omega|\bar\rho_t|\,|\eta|\,|\eta_t|\,dx\,ds
\le\gamma\int_0^t\|\eta_t(s)\|_H^2\,ds+\frac{C_6}{\gamma}\int_0^t
\|\bar\rho_t(s)\|^2_{L^4(\Omega)}\|\eta(s)\|^2_{L^4(\Omega)}\,ds
\nonumber\\
&\hspace*{-6mm}&\quad\le \gamma\int_0^t\|\eta_t(s)\|_H^2\,ds\,+\,\frac{C_7}{\gamma}\int_0^t\,
\|\bar\rho_t(s)\|^2_V\,\|\eta(s)\|^2_V\,ds \,.
\end{eqnarray}
Observe that by (\ref{eq:2.6}) the mapping $\,s\mapsto \|\bar\rho_t(s)\|^2_V\,$ belongs to $L^1(0,T)$.

Finally, we have, owing to the continuity of the embedding $W\subset L^\infty(\Omega)$ and
(\ref{eq:3.44}), 
\begin{eqnarray}
\label{eq:3.49}
&&\int_0^t\int_\Omega
2\,|\bar\mu_t|\,|\xi|\,|\eta_t|\,dx\,ds\,\nonumber\\[1mm]
&&\le\,\gamma\int_0^t\|\eta_t(s)\|_H^2\,ds\,+\,
\,\frac{C_8}{\gamma}\,\int_0^t\,
\|\bar\mu_t(s)\|^2_H\,\|\xi(s)\|_{L^\infty(\Omega)}^2ds\nonumber\\[1mm]
&&\le\,\gamma\int_0^t\!\!\|\eta_t(s)\|_H^2\,ds\,+\,\frac{C_9}{\gamma}\Bigl(\int_0^T\!\!\|h(s)\|^2_H\,ds\,+\,
\int_0^t\!\!\|\bar\mu_t(s)\|_H^2\,\|\Delta\xi(s)\|_H^2\,
ds\Bigr)\,,\nonumber\\
&&
\end{eqnarray}
where, owing to (\ref{eq:2.6}), the mapping $\,s\mapsto\|\bar\mu_t(s)\|_H^2\,$ belongs to $L^1(0,T)$.

Next, we test (\ref{eq:3.31}) by $\,-\Delta\xi_t\,$ to obtain, for every $t\in [0,T]$,
\begin{eqnarray}
\label{eq:3.50}
\delta\!\!\int_0^t\!\|\nabla\xi_t(s)\|_H^2\,ds+\frac 1 2 \|\Delta\xi(t)\|_H^2 =\int_0^t\!\!\int_\Omega
\left(-\left({\cal T}_\tau\eta\right)+f''(\bar\rho)\,\xi\right)\,
\Delta\xi_t\,dx\,ds.\quad
\end{eqnarray}

Now, by virtue of (\ref{eq:3.44}) and (\ref{eq:3.45}), and invoking Young's inequality, we have
\begin{eqnarray}
\label{eq:3.51}
&&\Bigl|\int_0^t\!\!\int_\Omega \left({\cal T}_\tau \eta\right) 
\Delta\xi_t\,dx\,ds\Bigr|\nonumber\\[1mm]
&&\le\,
\int_\Omega\!\! \left|\left({\cal T}_\tau \eta\right)(t)\right|\,|\Delta\xi(t)|\,dx\,
+\int_0^t\!\!\int_\Omega \left|\partial_t\left({\cal T}_\tau\eta\right)\right|\,|\Delta\xi|\,dx\,ds
\nonumber\\[1mm]
&&\le\,\frac 1 8\,\|\Delta\xi(t)\|_H^2\,+\,\gamma\int_0^t\!\|\eta_t(s)\|_H^2\,ds\,+\,C_{10}\,\Bigl(1+\frac 1
{\gamma}\Bigr)\int_0^T\!\!\|h(s)\|^2_H\,ds.\qquad\quad
\end{eqnarray}
Moreover, it turns out that
\begin{eqnarray}
\label{eq:3.52}
\hspace*{-8mm}\Bigl|\int_0^t\!\!\int_\Omega f''(\bar\rho)\,\xi\,\Delta\xi_t\,dx\,ds\Bigr|&\!\!\le\!\!&
\int_\Omega\!\!|f''(\bar\rho(t))|\,|\xi(t)|\,|\Delta\xi(t)|\,dx\nonumber\\[1mm]
&&+\int_0^t\!\!\int_\Omega \left|f'''(\bar\rho)\,\bar\rho_t\,\xi + f''(\bar\rho)\,\xi_t\right|
\,|\Delta\xi|\,dx\,ds.\quad
\end{eqnarray}

We have, owing to (\ref{eq:2.5}) and (\ref{eq:3.44}),
\begin{equation}
\label{eq:3.53}
\int_\Omega\!\!|f''(\bar\rho(t))|\,|\xi(t)|\,|\Delta\xi(t)|\,dx\,\le\,\frac 1 8 \,\|\Delta\xi(t)\|_H^2
\,+\,C_{11}\int_0^T\!\!\|h(s)\|^2_H\,ds\,.
\end{equation}

Also the second integral on the right-hand side of (\ref{eq:3.52}) is bounded, since (\ref{eq:2.5}), (\ref{eq:2.6}), (\ref{eq:3.44}), and (\ref{eq:3.45}) imply that 
\begin{eqnarray}
\label{pier:3}
\hskip-1.5cm &&\int_0^t\!\!\int_\Omega \left|f'''(\bar\rho)\,\bar\rho_t\,\xi + f''(\bar\rho)\,\xi_t\right| \,|\Delta\xi|\,dx\,ds \nonumber \\[1mm]
\hskip-1.5cm && \leq \, 
C_{12} \,\int_0^t\!\! \|\bar\rho_t (s)\|_{L^4(\Omega)}^2 \, \|\xi(s)\|_{L^4(\Omega)}^2 ds \, + \, \int_0^t\!\! \|\Delta\xi (s)\|_{H}^2 \,ds\, \nonumber \\[1mm]
\hskip-1.5cm && \leq \, 
C_{13} \max_{0\le t\le T}\,\|\xi(t)\|_{V}^2 \int_0^t\!\! \|\bar\rho_t (s)\|_{V}^2 \,ds \, + \, \int_0^t\!\! \|\Delta\xi (s)\|_{H}^2 \,ds\, 
\nonumber \\[1mm]
\hskip-2cm &&
\leq \, C_{14} \int_0^T\|h(s)\|_H^2\,ds\,,
\end{eqnarray}
thanks to the continuity of the embedding $V\subset L^4(\Omega)$.
Thus, combining the estimates (\ref{eq:3.46}) --(\ref{pier:3}), choosing $\gamma>0$ sufficiently
small, and invoking Gronwall's inequality, we can infer that
\begin{eqnarray}
\label{eq:3.54}
&{\displaystyle
  \int_0^T\!\!\Bigl(\|\eta_t(t)\|_H^2\,+\,\|\xi_t(t)\|_V^2\Bigr)
\,dt\,+\,
\max_{0\le t\le
  T}\,\left(\|\eta(t)\|_V^2\,+\,\|\xi(t)\|_W^2\right)}&\nonumber\\[1mm]
&{\displaystyle\le\,C_{15}\int_0^T\|h(t)\|_H^2\,dt\,.}&
\end{eqnarray}

Now, testing (\ref{eq:3.30}) by $\,-\Delta\eta$, and arguing as for (\ref{eq:3.46})--(\ref{eq:3.49}), we 
find that
$$ \int_0^T \|\eta(t)\|_W^2\,dt\,\le\,C_{16}\int_0^T\|h(t)\|_H^2\,dt\,.$$

Next, we differentiate Eq. (\ref{eq:3.31}) with respect to $t$. We obtain:
\begin{equation}
\label{eq:3.55}
\delta\,\xi_{tt}-\Delta\xi_t=\partial_t({\cal
  T}_\tau\eta)-f'''(\bar\rho)\,
\bar\rho_t\,\xi
-f''(\bar\rho)\,\xi_t\,\quad\mbox{a.\,e. in }\,Q.
\end{equation}

From (\ref{eq:2.5}), (\ref{eq:2.6}), (\ref{eq:3.44}), and (\ref{eq:3.54}), we can infer that the expression on the
right-hand side of (\ref{eq:3.55}) is bounded in $L^2(Q)$. Therefore, we may test (\ref{eq:3.55})
by any of the functions $\,\xi_t$, $-\Delta\xi_t$, and $\,\xi_{tt}$, in order to find that
\begin{equation}
\label{eq:3.56}
\int_0^T\!\!\Bigl(\|\xi_{tt}(t)\|_H^2\,+\,\|\Delta\xi_t(t)\|_H^2\Bigr)\,dt\,+\,
\max_{0\le t\le T}\,\|\xi_t(t)\|_V^2\,
\le\,C_{17}\int_0^T\|h(t)\|_H^2\,dt\,.
\end{equation}  

\vspace{2mm}
\underline{Step 3: \,\,Passage to the limit.} \,\quad Let $(\xi_N,\eta_N)$ denote the solution
to the system (\ref{eq:3.30})--(\ref{eq:3.33}) associated with $\tau_N=T/N$, for $N\in\nz$. In Step 2, we have shown that
there is some $C>0$, which does not depend on $N$, such that
\begin{eqnarray}
\label{eq:3.57}
&\|\xi_N\|_{H^2(0,T;H)\cap C^1([0,T];V)\cap H^1(0,T;W)\cap C(\overline
  Q)}&\nonumber\\[1mm]
&+\,
\|\eta_N\|_{H^1(0,T;H)\cap C^0([0,T];V)\cap L^2(0,T;W)}\,\le\,C.&
\end{eqnarray}
Hence, there is a subsequence, which is again indexed by $N$, such that 
\begin{eqnarray}
\label{eq:3.58}
\hspace*{-5mm}\xi_N\to \xi&\hspace*{-2mm}& \mbox{weakly star in }\,\,\,
H^2(0,T;H)\cap W^{1,\infty}(0,T;V)\cap H^1(0,T;W),\nonumber\\[2mm]
\hspace*{-5mm}\eta_N\to \eta&\hspace*{-2mm}&\mbox{weakly star in }\,\,\, H^1(0,T;H)\cap L^\infty(0,T;V)\cap L^2(0,T;W).
\end{eqnarray}
By compact embedding, we also have, in particular,
\begin{equation}
\label{eq:3.59}
\xi_N\to \xi \,\quad\mbox{strongly in }\,C(\overline{Q}), \qquad \eta_N\to\eta\quad\,\mbox{strongly in }\,
L^2(Q),
\end{equation}
so that $\,\bar\rho\,\eta_{N,t}\to\bar\rho\,\eta_t\,$ and $\,\bar\mu\,\xi_{N,t}\to\bar\mu\,\xi_t$,
both weakly in $\,L^2(Q)$, $\,f''(\bar\rho)\,\xi_N\to f''(\bar\rho)\,\xi\,$ strongly in $\,L^2(Q)$,
as well as $\,\bar\mu_t\,\xi_{N,t}\to \bar\mu_t\,\xi_t\,$ and $\,\bar\rho_t\,\eta_N\to
\bar\rho_t\,\eta$, both strongly in $L^1(Q)$. Finally, it is easily verified that $\,\{{\cal T}_{T/N}
\eta_N\}\,$ converges strongly in $L^2(Q)$ to $\eta$. In conclusion, we may pass to the limit as
$N\to\infty$ in the system (\ref{eq:3.30})--(\ref{eq:3.33}) (written for $\tau=T/N$) to find that the pair $(\xi,\eta)$
is in fact a strong solution to the linearized system (\ref{eq:3.21})--(\ref{eq:3.24}).

 It remains to show the uniqueness. But if $(\xi_1,
\eta_1)$, $(\xi_2,\eta_2)$ are two solutions having the above properties, then the pair $(\xi,\eta)$, where
$\xi=\xi_1-\xi_2$ and $\eta=\eta_1-\eta_2$, satisfies (\ref{eq:3.21})--(\ref{eq:3.24}) with $h= 0$. We thus may repeat
the first a priori estimate in Step 2 to conclude that $\xi=\eta=0$. This concludes the proof. \qed

\subsection{Directional differentiability of the control-to-state mapping}

In this section, we prove the following result.

\vspace{3mm}
{\bf Proposition 3.3} \,\quad {\em Suppose that the assumptions} (A1)--(A3) {\em are satisfied.
Then the solution operator $S$, viewed
as a mapping from $U_{ad}$, subset of $L^2(Q)$, into 
$$\,\left(H^1(0,T;H)\cap C^0([0,T];V)\cap L^2(0,T;W)\right)
\times \left(C^0([0,T];H)\cap L^2(0,T;V)\right),$$  is
directionally differentiable at $\bar u$ in the direction $h$. The directional derivative
$(\xi,\eta)=DS(\bar u)h$ is given by the unique solution $(\xi,\eta)$ to the linearized system}
(\ref{eq:3.21})--(\ref{eq:3.24}).

\vspace{3mm}
{\em Proof.} \,\quad Let $\bar\lambda>0$ be such that $\bar u+\lambda h\in U_{ad}$ for $0<\lambda\le\bar\lambda$. We put 
$$u^\lambda=\bar u+\lambda h,\,\quad (\rho^\lambda,\mu^\lambda)=S(u^\lambda),\,\quad
y^\lambda=\rho^\lambda-\bar\rho-\lambda\xi,\,\quad z^\lambda=\mu^\lambda-\bar\mu-\lambda\eta.$$ 
We have to show that there is a function 
$Z:[0,\bar\lambda]\to [0,+\infty)$ with\linebreak $\,\lim_{\lambda\searrow 0}\,Z(\lambda)/\lambda^2=0\,$
such that
\begin{equation}
\label{eq:3.60}
\|y^\lambda\|^2_{H^1(0,T;H)\cap C^0([0,T];V)\cap L^2(0,T;W)}\,+\,
\|z^\lambda\|^2_{C^0([0,T];H)\cap L^2(0,T;V)}\,\le \,Z(\lambda).
\end{equation}
Using the state system (\ref{eq:1.2})--(\ref{eq:1.5}) and the linearized system 
(\ref{eq:3.21})--(\ref{eq:3.24}), we easily verify
that for $0<\lambda\le\bar\lambda$ the pair $(y^\lambda,z^\lambda)$ is a strong solution
to the system
\begin{eqnarray}
\label{eq:3.61}
(\ve+2\bar\rho)\,z^\lambda_t+\bar\rho_t\,z^\lambda+\bar\mu\,y^\lambda_t
+2\bar\mu_t\,y^\lambda-\Delta z^\lambda\hspace{3.5cm}\nonumber\\
=
-2\left(\mu^\lambda_t-\bar\mu_t\right)\left(\rho^\lambda-\bar\rho\right)-
\left(\rho^\lambda_t-\bar\rho_t\right)\left(\mu^\lambda-\bar\mu\right)\, \quad
\mbox{a.\,e. in }\,Q,\\[1mm]
\label{eq:3.62}
\delta y^\lambda_t-\Delta y^\lambda+f'(\rho^\lambda)-f'(\bar\rho)-\lambda\, f''(\bar\rho)\,\xi
=z^\lambda,\,\quad\mbox{a.\,e. in }\,Q,\\[1mm]
\label{eq:3.63}
\frac{\partial y^\lambda}{\partial {\bf n}} =\frac{\partial z^\lambda}{\partial {\bf n}}=0,
\,\quad\mbox{a.\,e. on }\,\Sigma,\\[1mm]
\label{eq:3.64}
y^\lambda(x,0)=z^\lambda(x,0)=0 \,\quad\mbox{a.\,e. in }\,\Omega .
\end{eqnarray}
Notice that 
\begin{eqnarray*}
&&y^\lambda\in H^1(0,T;H)\cap C^0([0,T];V)\cap L^2(0,T;W)\cap C(\bar Q),\\[1mm]
&&z^\lambda\in H^1(0,T;H)\cap C^0([0,T];V)\cap L^2(0,T;W).
\end{eqnarray*}
For the sake of a better readability, in the following estimates we omit the superscript $\lambda$ of $y^\lambda$ and
$z^\lambda$. 
 As before, we denote by $C_i$ ($i\in\nz$) certain positive
constants that only depend
on  $\,\varepsilon, \delta, \rho_*, \rho^*,\mu^*, T, K^*_1, K^*_2, K^*_3$, but not on $\lambda$.

We now add $y$ on both sides of Eq. (\ref{eq:3.62}) and test the resulting equation by $y_t$. Using Young's
inequality, we find that for all $t\in [0,T]$ it holds
\begin{eqnarray}
\label{eq:3.65}
\frac {\delta} {2}\int_0^t \|y_t(s)\|_H^2\,ds +\frac 1 2 \left(\|\nabla y(t)\|_H^2+\|y(t)\|_H^2\right)
\,\le\,\frac{2}{\delta}\int_0^t \|z(s)\|_H^2\,ds \qquad\nonumber\\
\,+\,C_1 \int_0^t\|y(s)\|^2_H\,ds\,+\,C_2\int_0^t\|(f'(\rho^\lambda)-f'(\bar\rho)-\lambda f''(\bar\rho)
\,\xi)(s)\|_H^2\,ds\,.
\end{eqnarray} 
In order to handle the third term on the right-hand side of (\ref{eq:3.65}), we note that the stability 
estimate (\ref{eq:3.1}) implies, in particular, that
\begin{equation}
\label{eq:3.66}
\|\rho^\lambda-\bar\rho\|_{L^\infty(Q)}^2\,\le\,K^*_3\,\lambda^2\,\|h\|^2_{L^2(Q)}\,,
\end{equation}
that is, $\rho^\lambda\to\bar\rho$ uniformly on $\overline{Q}$ as $\lambda\searrow 0$. Since $f\in C^3(0,1)$,
we can infer from Taylor's theorem that
\begin{equation}
\label{eq:3.67}
\left|f'(\rho^\lambda)-f'(\bar\rho)-\lambda\,f''(\bar\rho)\,\xi\right|\le
\frac 1
2\,\max_{\rho_*\le\sigma\le\rho^*}\,\left|f'''(\sigma)\right|\,
\left|\rho^\lambda-\bar\rho\right|^2
+|f''(\bar\rho)|\,|y|\,\,\,\,\mbox{on }\,\overline Q.
\end{equation}
 It then follows from the estimates (\ref{eq:3.1}) and  (\ref{eq:3.65}) that 
\begin{eqnarray}
\label{eq:3.68}
\hspace*{-4mm}&&{\displaystyle \frac {\delta} {2}\int_0^t\!\!\! \|y_t(s)\|_H^2\,ds 
+\frac 1 2 \|y(t)\|_V^2
\,\le\,\frac{2}{\delta}\int_0^t\!\!\! \|z(s)\|_H^2\,ds 
\,+\,C_3\!\! \int_0^t\!\!\!\|y(s)\|^2_H\,ds\,+\,C_4\lambda^4,}\nonumber\\
\hspace*{-4mm}&&
\end{eqnarray}

since we can assume that $\,\|h\|_{L^2(Q)}\le 1$.
Next, observe that
$\,2\,\bar\rho\,z\,z_t=\left(\bar\rho\,z^2\right)_t-\bar\rho_t\,z^2$.
Therefore, 
testing
(\ref{eq:3.61}) by $z$ yields for every $t\in [0,T]$ that
\begin{eqnarray}
\label{eq:3.69}
&\hspace*{-8mm}&{\displaystyle \int_\Omega\left(\frac {\ve} 2 +\bar\rho(t)\right)z^2(t)\,dx\,+\,\int_0^t\|\nabla z(s)\|_H^2\,ds
\,=\,-\int_0^t\!\!\int_\Omega\left(\bar\mu\,y_t\,+\,2\,\bar\mu_t
\,y\right)z\,dx\,ds}\nonumber\\
&\hspace*{-8mm}&{\displaystyle -2\int_0^t\!\!\int_\Omega \left(\mu^\lambda_t-
\bar\mu_t\right)\left(\rho^\lambda-\bar\rho\right)\,z
\,dx\,ds -\int_0^t\!\!\int_\Omega 
\left(\rho^\lambda_t-\bar\rho_t\right)\left(\mu^\lambda-\bar\mu\right)
\,z\,dx\,ds.}\nonumber\\
&\hspace*{-8mm}&
\end{eqnarray} 
We estimate the terms on the right-hand side of (\ref{eq:3.69}) individually. At first, 
using (\ref{eq:2.6}) and Young's
inequality, we find that
\begin{equation}
\label{eq:3.70}
\int_0^t\!\!\int_\Omega |\bar\mu|\,|y_t|\,|z|\,dx\,ds\,\le\,\gamma\int_0^t\|y_t(s)\|_H^2\,ds
\,+\,\frac{C_5}{\gamma}\int_0^t\|z(s)\|_H^2\,ds.
\end{equation}
Moreover, using the continuity of the embedding $H^1(\Omega)\subset L^4(\Omega)$, 
as well as H\"older's and Young's inequalities,
\begin{eqnarray}
\label{eq:3.71}
 \int_0^t\!\!\int_\Omega |\bar\mu_t|\,|y|\,|z|\,dx\,ds\,\le\,\int_0^t\|\bar\mu_t(s)\|_H
 \,\|z(s)\|_{L^4(\Omega)}\,\|y(s)\|_{L^4(\Omega)}\,ds\nonumber\\
 \le \,\gamma\int_0^t\|z(s)\|_V^2\,ds\,+\,\frac{C_6}{\gamma}\int_0^t\|\bar\mu_t(s)\|^2_H\,\|
 y(s)\|_V^2\,ds\,.
 \end{eqnarray}
 Observe that by (\ref{eq:2.6}) the mapping $\,s\mapsto \|\bar\mu_t(s)\|_H^2\,$ belongs to $L^1(0,T)$.
 
At this point, we can conclude from (\ref{eq:3.1}) and (\ref{eq:3.66}), invoking Young's inequality, that
 \begin{eqnarray}
 \label{eq:3.72}
 &&\int_0^t\!\!\int_\Omega 2\,\left|\mu_t^\lambda-\bar\mu_t\right|\,\left|\rho^\lambda
 -\bar\rho \right|\,|z|\,dx\,ds\nonumber\\
 &&\le\,2\int_0^t \left\|(\mu_t^\lambda-\bar\mu_t)(s)\right\|_H\,
 \left\|(\rho^\lambda  -\bar\rho)(s) \right\|_{L^\infty(Q)}\,\|z(s)\|_H\,ds\nonumber\\
&& \le\,C_7\,\left\|(\rho^\lambda  -\bar\rho)(s) \right\|_{L^\infty(Q)}^2\int_0^t
 \left\|(\mu_t^\lambda-\bar\mu_t)(s)\right\|_H^2\,ds \,+\,\int_0^t\|z(s)\|^2_H\,ds\nonumber\\
 &&\le \int_0^t\|z(s)\|^2_H\,ds\,+\,C_8\,\lambda^4\,.
 \end{eqnarray}
 
 Finally, we invoke (\ref{eq:3.1}) and H\"older's and Young's inequalities, as well as the continuity of the 
 embedding $H^1(\Omega) \subset L^4(\Omega)$, to obtain that
  \begin{eqnarray}
  \label{eq:3.73}
 &&\int_0^t\!\!\int_\Omega \left|\rho_t^\lambda-\bar\rho_t\right|\,\left|\mu^\lambda
 -\bar\mu \right|\,|z|\,dx\,ds\nonumber\\
 &&\le \,\max_{0\le s\le t}\|z(s)\|_H \int_0^t\left\|(\rho_t^\lambda-\bar\rho_t)(s)\right\|_{L^4(\Omega)}\,\left\|(\mu^\lambda
 -\bar\mu)(s) \right\|_{L^4(\Omega)}\,ds\nonumber\\
 &&\le\,\gamma \,\max_{0\le s\le t}\|z(s)\|_H^2\,+\,\frac{C_9}{\gamma}\int_0^t\left \|
 (\rho_t^\lambda-\bar\rho_t)(s)\right\|_V^2\,ds\,\,
 \int_0^t\left \|
 (\mu^\lambda-\bar\mu)(s)\right\|_V^2\,ds\nonumber\\[2mm]
 &&\le \,\gamma \,\max_{0\le s\le t}\|z(s)\|_H^2\,+\,C_{10}\,\lambda^4\,.
 \end{eqnarray}
 
 Combining the estimates (\ref{eq:3.68})--(\ref{eq:3.73}), taking the maximum with respect to $t\in [0,T]$,
 adjusting $\gamma>0$ appropriately small, and invoking Gronwall's lemma, we arrive at the conclusion that $(y^\lambda,z^\lambda)
 =(y,z)$ satisfies the inequality
 \begin{equation}
 \label{eq:3.74}
 \|y^\lambda\|_{H^1(0,T;H)\cap C^0([0,T];V)}^2\,+\,\|z^\lambda\|^2_{C^0([0,T];H)\cap L^2([0,T];V)}
 \,\le\,C_{11}\,\lambda^4\,.
 \end{equation}
 Finally, testing (\ref{eq:3.62}) by $\,-\Delta y^\lambda$, and using (\ref{eq:3.67}), we find that also
 \begin{equation}
 \label{eq:3.75}
 \|y^\lambda\|^2_{L^2(0,T;W)}\,\le\, C_{12}\,\lambda^4\,.
 \end{equation}
 This concludes the proof of the assertion. \qed
 
 \vspace{7mm}
 {\bf Corollary 3.4} \,\quad {\em Let the assumptions} (A1)--(A3) {\em be fulfilled, and
 let $\bar u\in U_{ad}$  be an optimal control for the problem} {\bf (CP)} {\em with associated
 state $(\bar\rho,\bar\mu)=S(\bar u)$. Then, for every $v\in U_{ad}$, 
 \begin{equation}
 \label{eq:3.76}
 \int_0^T\!\!\!\int_ \Omega\!\! \beta_2\,\bar u(v-\bar u)\,dx\,dt +\int_\Omega\!\! (\bar\rho(T)-\rho_T)\,
 \xi(T)\,dx\,+\,\int_0^T\!\!\int_\Omega\!\!\beta_1\,(\bar\mu-\mu_T)\,\eta\,dx\,dt\,\ge\,0,\quad
  \end{equation} 
 where $(\xi,\eta)$ is the unique solution to the linearized system} 
 (\ref{eq:3.21})--(\ref{eq:3.24}) {\em associated with 
 $h=v-\bar u$.} 
 
 \vspace{2mm}
 {\em Proof.} \,\quad Let $\,v\in U_{ad}$ be arbitrary. Then $h=v-\bar u$ is an admissible
 direction, since $\bar u +\lambda h  \in U_{ad}$ for $0<\lambda\le
 1$. 
For any such $\lambda$, we have
 \begin{eqnarray*}
 &\hspace*{-6mm}&0\,\le \,\frac{J(\bar u+\lambda h,S(\bar u+\lambda h))-J(\bar u,S(\bar u))}\lambda\nonumber\\[3mm]
 &\hspace*{-6mm}&\le\, \frac{J(\bar u+\lambda h,S(\bar u+\lambda h))
-J(\bar u,S(\bar u+\lambda h))}{\lambda}+
 \frac{ J(\bar u,S(\bar u+\lambda h))-J(\bar u,S(\bar u))}{\lambda}\,.\qquad
 \end{eqnarray*}
 It follows immediately from the definition of the cost functional $J$ that the first summand on the 
 right-hand side of this inequality converges to $\,\int_0^T\!\!\int_\Omega \beta_2\,\bar u\,h\,dx\,dt\,$
 as $\lambda\searrow 0$. For the second summand, we obtain from Proposition 3.3 that
 \begin{eqnarray*}
&{\displaystyle \lim_{\lambda\searrow 0}\,\frac{ J(\bar u,S(\bar u+\lambda h))-J(\bar 
u,S(\bar u))}{\lambda}}&\\
&{\displaystyle  = \int_\Omega\!\! (\bar\rho(T)-\rho_T)\,
 \xi(T)\,dx\,+\,\int_0^T\!\!\int_\Omega\!\!\beta_1\,(\bar\mu-\mu_T)\,
\eta\,dx\,dt\,,}&
 \end{eqnarray*}
 whence the assertion follows. \qed
 
 \subsection{The optimality system}
 
 Let $\bar u\in U_{ad}$ be an optimal control for {\bf (CP)} with associated state 
 $(\bar\rho,\bar\mu)=S(\bar u)$. Then, for every $\,v\in U_{ad}$, (\ref{eq:3.76}) holds. We 
 now aim to eliminate $(\xi,\eta)$ by introducing the adjoint
 state variables. To this end, we consider the {\em adjoint system}\,:
 
 \begin{eqnarray}
 \label{eq:3.77}
 -(\ve+2\bar\rho)\,q_t-\bar\rho_t\,q-\Delta q=p+\beta_1\left(\bar\mu-\mu_T\right)\,\quad
 \mbox{a.\,e. in }\,Q,\\[1mm]
 \label{eq:3.78}
 \frac{\partial q}{\partial \bf n}=0\quad \mbox{a.\,e. in }\,\Sigma,\qquad q(x,T)=0\, \quad
 \mbox{a.\,e. in }\,\Omega,\\[2mm]
 \label{eq:3.79}
 -\delta p_t-\Delta p+f''(\bar\rho)\,p=\bar\mu\,q_t-\bar\mu_t\,q\,\quad
 \mbox{in }\,Q,\\[1mm]
 \label{eq:3.80}
 \frac{\partial p}{\partial \bf n}=0\quad \mbox{on }\,\Sigma,\qquad  
 \delta\,p(T)=\bar\rho(T)-\rho_T \,\quad \mbox{in }\,\Omega\,, 
 \end{eqnarray}
 
 which is a linear backward-in-time parabolic system for the adjoint state variables $p$ and $q$.
 
 It must be expected that the adjoint state variables $(p,q)$ be less regular than the state
 variables $(\bar\rho,\bar\mu)$. Indeed, we only have $p(T)\in L^2(\Omega)$, and thus 
 (\ref{eq:3.79}) and (\ref{eq:3.80})
 should be interpreted in the ususal weak sense. That is, we look for a vector-valued function
 $\,p\in H^1(0,T;V^*)\cap C^0([0,T];H)\cap L^2(0,T;V)$ that, in addition to the final time condition
 (\ref{eq:3.80}), satisfies 
 \begin{eqnarray}
 \label{eq:3.81}
 &&\langle -\delta\,p_t(t),v\rangle_{V^*,V}\,+\,\int_\Omega\nabla p(t)\cdot\nabla v\,dx \,+\,
 \int_\Omega f''(\bar\rho(t))\,p(t)\,v\,dx\nonumber\\[1mm]
 &&\quad= \int_\Omega \left(\bar\mu(t)\,q_t(t)-\bar\mu_t(t)\,q(t)\right)\,v\,dx\,,
 \end{eqnarray}
 for every $v\in V$ and almost every $t\in (0,T)$. 
  Notice that if $q\in H^1(0,T;H)\cap C^0([0,T];V)$, then it is easily seen that 
 $\,\bar\mu\,q_t-\bar\mu_t\,q\in L^{3/2}(Q)$, so that the integral on the right-hand side of
 (\ref{eq:3.81}) makes sense. On the other hand, if $p$ has the expected regularity then
 the solution to (\ref{eq:3.77}), (\ref{eq:3.78}) should belong to 
 $H^1(0,T;H)\cap C^0([0,T];V)\cap L^2(0,T;W)$.
 
 \vspace{5mm}
 {\bf Lemma 3.5} \,\quad {\em Suppose that the system} (\ref{eq:3.77})--(\ref{eq:3.80}) 
 {\em has a unique solution
 $(p,q)$ where $p\in H^1(0,T; V^*)\cap C^0([0,T];H)\cap L^2(0,T;V)$ and $q\in
 H^1(0,T;H)\cap C^0([0,T];V)\cap L^2(0,T;W)$. Then we have}
  \begin{equation}
  \label{eq:3.82}
 \int_\Omega (\bar\rho(x,T)-\rho_T(x))\,
 \xi(x,T)\,dx\,+\,\int_0^T\!\!\int_\Omega \beta_1\,(\bar\mu-\mu_T)\,\eta\,dx\,dt = \int_0^T\!\!
 \int_\Omega q\,h\,dx\,dt\,.
 \end{equation} 
 
 \vspace{2mm}
 {\em Proof.} \,\quad The assertion follows from repeated integration by parts, using the well-known integration by parts formula
 $$\int_0^T\bigl(\langle v_t(t),w(t)\rangle_{V^*,V}\,+\, \langle w_t(t),v(t)\rangle_{V^*,V}\bigr)\,dt
 =\int_\Omega\bigl(v(T)w(T)-v(0)w(0)\bigr)\,dx,$$
 which holds for all functions $v,w\in H^1(0,T;V^*)\cap L^2(0,T;V)$. Since this calculation
 is standard in optimal control theory, we may leave it to the reader to work out the details.\qed
 
 \vspace{7mm}
 {\bf Proposition 3.6} \,\quad {\em The adjoint system} (\ref{eq:3.77})--(\ref{eq:3.80}) {\em has a unique 
 solution $(p,q)$ with $p\in H^1(0,T; V^*)\cap C^0([0,T];H)\cap L^2(0,T;V)$ and $q\in
 H^1(0,T;H)\cap C^0([0,T];V)\cap L^2(0,T;W)$, where} (\ref{eq:3.79}) and (\ref{eq:3.80})$_1$ 
 {\em are understood
 in the sense of} (\ref{eq:3.81}). 
  
 \vspace{2mm}
 {\em Proof.} \,\quad We proceed in a series of steps.
 
 \vspace{2mm}
 \underline{Step 1: \,Approximation.} \,\quad As in the proof of Proposition 3.2, 
 we employ a delay technique. However, this  time we have to use a {\em negative} delay, 
 since the system (\ref{eq:3.77})--(\ref{eq:3.80}) runs backwards in time.
 
 For $0<\tau<T$, we define the translation operator $\tilde{\cal T}_\tau:L^1(0,T;H)\to
 L^1(0,T;H)$ by putting, for every $v\in L^1(0,T;H)$ and almost all $t\in (0,T)$,

\vspace{-1mm} 
\begin{equation}
 \label{eq:3.83}
 \bigl(\tilde{\cal T}_\tau v\bigr)(t)=v(t+\tau)\,\quad\mbox{if }\,t\le T-\tau,\quad\mbox{ and }\quad
 \bigl(\tilde{\cal T}_\tau v\bigr)(t)=0\,\quad\mbox{if }\,t> T-\tau;
 \end{equation}

\vspace{3mm} 
 clearly, we have that $\|\tilde{\cal T}_\tau v\|_{L^2(Q)}\le\|v\|_{L^2(Q)}$,
  for every $v\in L^2(Q)$ and any $\tau\in (0,T)$.

We then consider for $\tau\in (0,T)$ the following approximating problem: find functions 
\begin{eqnarray}
\label{eq:3.84}
&&p^\tau\in H^1(0,T;V^*)\cap C^0([0,T];H)\cap L^2(0,T;V),\nonumber\\
&&q^\tau\in H^1(0,T;H)\cap
C^0([0,T];V)\cap L^2(0,T;W)
\end{eqnarray}
that solve the system 
\begin{eqnarray}
\label{eq:3.85}
&-(\ve+2\bar\rho)\,q^\tau_t-\Delta q^\tau + q^\tau = (1+\bar \rho_t)\,
\tilde{\cal T}_\tau q^\tau 
+\tilde{\cal T}_\tau p^\tau +\beta_1(\bar\mu-\mu_T)\,&\nonumber\\
&\mbox{a.\,e. in }\,Q,&\\[1mm]
\label{eq:3.86}
&{\displaystyle \frac{\partial q^\tau}{\partial \bf n}=0\,\quad
\mbox{a.\,e. on }\,\Sigma, \,\quad q^\tau(x,T)=0 \,\,\quad 
\mbox{a.\,e. in }\,\Omega,}&\\[1mm]
\label{eq:3.87}
&{\displaystyle \langle -\delta\,p^\tau_t(t),v\rangle_{V^*,V}+\int_\Omega\!\! 
\nabla p^\tau(t) \cdot\nabla v\,dx 
\,+\,\int_\Omega\!\!f''(\bar\rho(t))\,p^\tau(t)\,v\,dx}&\nonumber\\[1mm]
&{\displaystyle =\int_\Omega\!\!\Bigl(\bar\mu(t)\,q_t^\tau(t)-\bar\mu_t(t)q^\tau(t)\Bigr)
v\,dx \quad\forall\,v\in V,\quad\mbox{for a.\,e. }t\in (0,T),}\quad&\\[1mm]
\label{eq:3.88}
&{\displaystyle \delta\,p^\tau(T)=\bar\rho(T)-\rho_T\,\quad \mbox{a.\,e. in }\,
\Omega.}&
\end{eqnarray}
We choose for $\tau\in (0,T)$ the discrete values $\tau=T/N$, where $N\in\nz$ is arbitrary, and we put
$t_n=n\,\tau$, $0\le n\le N$, and  $I_n=(t_n,T)$. For $n=N-1,\ldots,1,0$, we solve the problem:
\begin{eqnarray}
\label{eq:3.89}
&-(\ve+2\bar\rho)\,q_{n,t}-\Delta q_n + q_n = (1+\bar \rho_t)\,\tilde{\cal T}_\tau q_n 
+\tilde{\cal T}_\tau p_n +\beta_1(\bar\mu-\mu_T)&\nonumber\\
&\qquad\qquad\mbox{a.\,e. in }\,\Omega\times I_n,&\\[1mm]
\label{eq:3.90}
&{\displaystyle \frac{\partial q_n}{\partial \bf n}=0\,\quad
\mbox{a.\,e. on }\,\Sigma, \,\quad q_n(x,T)=0 \,\,\quad
\mbox{a.\,e. in }\,
\Omega,}& \\[1mm]
\label{eq:3.91}
&{\displaystyle \langle -\delta\,p_{n,t}(t),v\rangle_{V^*,V}+\int_\Omega\!\! \nabla p_n(t) \cdot\nabla v\,dx 
\,+\,\int_\Omega\!\!f''(\bar\rho(t))\,p_n(t)\,v\,dx
}&\nonumber\\[1mm]
&{\displaystyle =\int_\Omega\!\!\Bigl(\bar\mu(t)\,q_{n,t}(t)-\bar\mu_t(t)\, 
q_n(t)\Bigr)\,v\,dx\,\quad\mbox{for all }\,v\in V,
\quad\mbox{for a.\,e. }\,t\in (t_n,T),}&\nonumber\\
&{}&\\
\label{eq:3.92}
&\delta\,p_n(T)=\bar\rho(T)-\rho_T\,\quad \mbox{a.\,e. in }\,\Omega.&
\end{eqnarray}

Here, $\tilde{\cal T}_\tau$ acts on functions that are not defined on the entire interval $(0,T)$; however,
for $n<N-1$ it is still defined by (\ref{eq:3.83}), while for $n=N-1$ we simply put $\tilde{\cal T}_\tau p_n
=\tilde{\cal T}_\tau q_n=0$.
Whenever the pairs $(p_k,q_k)$ with
\begin{eqnarray}
\label{eq:3.93}
p_k\in H^1(I_k;V^*)\cap C^0(\bar I_k;H)\cap L^2(I_k;V),\\
\label{eq:3.94}
q_k\in  H^1(I_k;H)\cap C^0(\bar I_k;V)\cap L^2(I_k;W),
\end{eqnarray}
 have been constructed for $0< n\le k\le N-1$, we look for the pair $(\xi_{n-1},\eta_{n-1})$ 
that coincides with $(\xi_n,\eta_n)$ in $I_n$. Note that 
 the linear parabolic problem (\ref{eq:3.89}), (\ref{eq:3.90}) has a 
unique solution $q_{n-1}$ on $\Omega\times I_{n-1}$ that satisfies (\ref{eq:3.94}) for $k=n-1$
(see, e.\,g., \cite{LSU}) . On inserting
$q_{n-1}$ in (\ref{eq:3.91}) (with $n$ replaced by $n-1$), we then find (e.\,g., by using an appropriate
Galerkin approximation) that the linear parabolic problem (\ref{eq:3.91}), (\ref{eq:3.92}) admits a
unique solution $p_{n-1}$ that fulfills (\ref{eq:3.93}) for $k=n-1$. Hence, we conclude that 
$(p^\tau ,q^\tau)=(p_0,q_0)$ satisfies (\ref{eq:3.85})--(\ref{eq:3.88}) and (\ref{eq:3.84}).  

\vspace{3mm}
\underline{Step 2: \,\,A priori estimates.} \,\quad We now prove a series of a priori estimates for the functions $(p^\tau,q^\tau)$. In the following, we denote by $C_i$ ($i\in\nz$) some generic positive constants, 
which may depend on $\,\varepsilon,\delta,\rho_*, \rho^*,\mu^*, T, K^*_1, K^*_2$, 
but not on $\tau$ (i.\,e., not on $N$). For the sake of simplicity,
we omit the superscript $\tau$ and simply write $(p,q)$. 

\vspace{2mm}
We multiply (\ref{eq:3.85}) by $-q_t$ and integrate over $\Omega\times [t,T]$ to obtain, using Young's 
inequality,
\begin{eqnarray}
\label{eq:3.95}
&\hspace*{-6mm}&\frac {\ve} 2\int_t^T\|q_t(s)\|_H^2\,ds\,+\,\frac 1 2 \left(\|\nabla q(t)\|_H^2\,+\,\|q(t)\|_H^2\right)
\nonumber\\
&\hspace*{-6mm}&\,\le C_1\,+\,C_2\int_t^T\left(\|p(s)\|_H^2\,+\,\|q(s)\|_H^2\right)ds\,+\,\int_t^T\!\!\int_\Omega
|\bar\rho_t|\,|\tilde{\cal T}_\tau q|\,|q_t|\,dx\,ds\,.\qquad\quad 
\end{eqnarray}
Moreover, by virtue of H\"older's and Young's inequalities, and invoking the continuity of the
embedding $H^1(\Omega)\subset L^4(\Omega)$, we have, for any $\gamma>0$, that
\begin{eqnarray}
\label{eq:3.96}
&&\int_t^T\!\!\int_\Omega
|\bar\rho_t|\,|\tilde{\cal T}_\tau q|\,|q_t|\,dx\,ds\,\le\,\int_t^T \|\bar\rho_t(s)\|_{L^4(\Omega)}\,\|\tilde{\cal T}_\tau q(s)\|_{L^4(\Omega)}\,\|q_t(s)\|_H\,ds\,\nonumber\\
&&\,\le\,\gamma\int_t^T\|q_t(s)\|_H^2\,ds\,+\,\frac{C_3}{\gamma}
\int_t^{T-\tau} \|\bar\rho_t(s)\|_V^2\,\|q(s+\tau)\|_V^2\,ds\,.
\end{eqnarray}
Observe that by (\ref{eq:2.6}) the mapping $\,s\mapsto \|\bar\rho_t(s)\|_V^2\,$ belongs to $L^1(0,T)$.

\vspace{2mm}
Next, we insert $v=p(t)$ in (\ref{eq:3.87}) and integrate over $[t,T]$, where $t\in [0,T]$. 
We find, using (\ref{eq:3.92}), that
\begin{eqnarray}
\label{eq:3.97}
&&\frac {\delta} 2 \,\|p(t)\|_H^2\,+\,\frac 1 2 \int_t^T\|\nabla p(s)\|_H^2\,ds\,\le\,C_4\,+\,C_5
\int_t^T\|p(s)\|_H^2\,ds\nonumber\\
&&\,+\,\int_t^T\!\!\int_\Omega |\bar\mu|\,|q_t|\,|p|\,dx\,ds\,+\,
\int_t^T\!\!\int_\Omega |\bar\mu_t|\,|q|\,|p|\,dx\,ds\,.
\end{eqnarray}
Let us denote for short the last two integrals on the right-hand side of (\ref{eq:3.97}) by $I_1$ and $I_2$,
 respectively. Since $\bar\mu\in L^\infty(Q)$, we have that
 \begin{equation}
 \label{eq:3.98}
 I_1\,\le\,\gamma\int_t^T\|q_t(s)\|_H^2\,ds\,+\,\frac{C_6}{\gamma}
\int_t^T \|p(s)\|_H^2\,ds\,.
 \end{equation}
 Moreover, we conclude from (\ref{eq:2.6}), invoking H\"older's and Young's inequalities, that
 \begin{eqnarray}
 \label{eq:3.99}
 &&I_2\,\le\,\int_t^T \|\bar\mu_t(s)\|_H\,\|q(s)\|_{L^4(\Omega)}\,\|p(s)\|_{L^4(\Omega)}\,ds\,\nonumber\\
 &&\,\le\,\gamma\int_t^T \|p(s)\|_V^2\,ds\,+\,\frac{C_7}{\gamma}
 \int_t^T \|\bar\mu_t(s)\|_H^2\,\|q(s)\|_V^2\,ds\,,
 \end{eqnarray}
 where the mapping $\,s\mapsto \|\bar\mu_t(s)\|_H^2\,$ belongs to $L^1(0,T)$.
 
 Now, we combine the estimates (\ref{eq:3.95})--(\ref{eq:3.99}). On choosing $\gamma>0$ sufficiently small,
 and on applying Gronwall's lemma, we find that
 \begin{equation}
 \label{eq:3.100}
 \|p\|_{C^0([0,T];H)\cap L^2(0,T;V)}\,+\,\|q\|_{H^1(0,T;H)\cap C^0([0,T];V)}\,\le \,C_8\,.
 \end{equation}
 It is now a standard matter to verify, by comparison in (\ref{eq:3.85}) and (\ref{eq:3.87}), respectively, 
  that also
 \begin{equation}
 \label{eq:3.101}
 \|p\|_{H^1(0,T;V^*)}\,+\,\|q\|_{L^2(0,T;W)}\,\le \,C_9\,.
 \end{equation}   
 
 \vspace{2mm}
\underline{Step 3: \,\,Passage to the limit.} \,\quad Let $(p_N,q_N)$ denote the solution
to the system (\ref{eq:3.85})--(\ref{eq:3.86}) associated with $\tau_N=T/N$, for $N\in\nz$. 
In Step 2, we have shown that
there is some $C>0$, which does not depend on $N$, such that
\begin{eqnarray*}
\|p_N\|_{H^1(0,T;V^*)\cap C^0([0,T];H)\cap L^2(0,T;V)}\,+\,
\|q_N\|_{H^1(0,T;H)\cap C^0([0,T];V)\cap L^2(0,T;W)}\,\le\,C.
\end{eqnarray*}
Hence, there is a subsequence, which is again indexed by $N$, such that 
\begin{eqnarray*}
&&p_N\to p,\,\quad\mbox{weakly star in }\,
H^1(0,T;V^*)\cap L^\infty(0,T;H)\cap L^2(0,T;V),\\[2mm]
&& q_N\to q,\,\quad\,\mbox{weakly star in }\,
H^1(0,T;H)\cap L^\infty(0,T;V)\cap L^2(0,T;W).
\end{eqnarray*}
By compact embedding, we also have the strong convergences (see, e.\,g., \cite{Simon}) 
\begin{eqnarray*}
&p_N\to p,\quad\mbox{strongly in }\,L^2(Q),&\\[1mm] 
&q_N\to q,\,\quad \mbox{strongly in }\,
C^0([0,T];H)\cap L^2(0,T;V).&
\end{eqnarray*}
From this, we can conclude the following convergences:
\begin{eqnarray*}
&&\bar\rho\,q_{N,t}\to \bar\rho \,q_t \quad\mbox{weakly in } L^2(Q), \,\quad \Delta q_N\to\Delta q
\quad\mbox{weakly in } L^2(Q),\nonumber\\[1mm]
&&\tilde{\cal T}_{T/N}q_N\to q \quad\mbox{strongly in } L^2(Q),\,\quad\tilde{\cal T}_{T/N}p_N\to p
\quad\mbox{strongly in } L^2(Q),\nonumber\\[1mm]
&&\bar\rho_t\,\tilde{\cal T}_{T/N}q_N\to \bar\rho_t\,q\quad\mbox{strongly in } L^1(Q).
\end{eqnarray*}
Hence, passing to the limit as $N\to\infty$ in (\ref{eq:3.85})--(\ref{eq:3.86}) for $\,\tau=T/N$, 
we find that the pair $(p,q)$ gives a strong solution to the parabolic problem 
(\ref{eq:3.77})--(\ref{eq:3.78}). Next, we notice that the
weak convergence of $\{p_N\}$ to $p$ in $H^1(0,T;V^*)\cap L^2(0,T;V)$ implies that
$\,p_N\to p\,$ weakly in $C^0([0,T];H)$. We may thus conclude, in particular, that
$\,\delta\,p(T)=\bar\rho(T)-\rho_T$. Since $f''(\bar\rho)$ and $\bar\mu$ are bounded, 
we also have the following convergences:
\begin{eqnarray*}
f''(\bar\rho)\,p_N\to f''(\bar\rho)\,p \quad\mbox{strongly in } L^2(Q),\,\quad
\bar\mu\,q_{N,t}\to \bar\mu\,q_t \quad\mbox{weakly in } L^2(Q).
\end{eqnarray*} 

Therefore, (\ref{eq:3.81}) is fulfilled for any $v$ in $C^1(\overline
\Omega)$, which is a dense subset of $H^1(\Omega)$, because $\Omega$ is a Lipschitz domain.
From this, it easily follows that (\ref{eq:3.81}) is satisfied for every $v\in H^1(\Omega)$. 
In conclusion, the pair $(p,q)$ is a solution to the
adjoint system (\ref{eq:3.77})--(\ref{eq:3.80}) that enjoys the asserted smoothness properties. 

 Uniqueness remains to be shown. But if $(p_1,
q_1)$, $(p_2,q_2)$ are two solutions having the above properties, then the pair $(p,q)$, where
$p=p_1-p_2$ and $q=q_1-q_2$, satisfies (\ref{eq:3.77})--(\ref{eq:3.80}), 
where the inhomogeneities $\beta_1(\bar\mu -\mu_T)$ and $\bar\rho-\rho_T$ on the right-hand sides 
of (\ref{eq:3.77}) and (\ref{eq:3.80}), respectively, cancel 
out by subtraction. 
We may then repeat the a priori estimates of Step 2 to see that in the present situation the constants
$C_1$ and $C_4$ appearing, respectively, in (\ref{eq:3.95}) and (\ref{eq:3.97}), simply do not occur. 
Consequently, the application of Gronwall's lemma yields 
$p=q=0$.
This concludes the proof.\\\hspace*{5cm}\hfill{\qed}

\vspace{4mm}
In summary, we have proved the following result concerning first-order 
necessary optimality conditions.

\vspace{4mm}
 {\bf Theorem 3.7} \,\quad {\em Suppose that $\bar u\in U_{ad}$ is an optimal control for} {\bf (CP)}
 {\em with associated state $(\bar\rho,\bar\mu)=S(\bar u)$. Then the adjoint system} 
 (\ref{eq:3.77})--(\ref{eq:3.80}) {\em
 has a unique weak solution $(p,q)$ with}
 $p\in H^1(0,T; V^*)\cap C^0([0,T];H)\cap L^2(0,T;V)\,,\quad q\in
 H^1(0,T;H)\cap C^0([0,T];V)\cap L^2(0,T;W)${\em ;
 moreover, for any $v\in U_{ad}$, we have the inequality:}
 \begin{equation}
 \label{eq:3.104}
 \int_0^T\!\!\int_\Omega \beta_2\,\bar u\,(v-\bar u)\,dx\,dt\,+\,\int_0^T \!\!\int_\Omega q\,(v-\bar u)
 \,dx\,dt\,\ge\,0\,.
 \end{equation} 
 
 \vspace{4mm}
 {\bf Remark:} \,\quad 5. Since $U_{ad}$ is a nonempty, closed and convex subset of $L^2(Q)$, 
 (\ref{eq:3.104}) has the following implications:
 \begin{itemize}
 \item For $\,\beta_2>0$, the optimal control $\bar u$ 
 is nothing but the $L^2(Q)$ orthogonal projection of $\,-\beta_2^{-1}\,q\,$ onto $U_{ad}$. In 
 other words, 
 $$\bar u(x,t)={\cal P}_{[0,U(x,t)]}\left(-\beta_2^{-1}\,q(x,t)\right)\,\quad\mbox{a.\,e. in }\,Q\,,
 $$
 where, for any  $\,a,b\in\rz\,$ such that $\,a\le b$,  
$$\,{\cal P}_{[a,b]}(u):=\mbox{min}\,\{b,\,
 \mbox{max}\,\{a,u\}\}.$$ 
 \item For $\beta_2=0$, we have that almost everywhere
 $$\bar u(x,t)=\left\{
 \begin{array}{ll}
 0&\mbox{if }\,q(x,t)>0\\[1mm]
 U(x,t)&\mbox{if }\,q(x,t)<0
  \end{array}
  \right. \,,
  $$
  a {\em bang-bang} situation where no information can be recovered for points at which 
  $q(x,t)=0$.
  \end{itemize}
\vspace{5mm}
\noindent
{\bf Acknowledgement. }
P. Colli and G. Gilardi 
gratefully acknowledge the financial support of
 the MIUR-PRIN Grant 2008ZKHAHN {\em ``Phase transitions, 
hysteresis and multiscaling''} and of the IMATI of CNR in Pavia. 
The work of J. Sprekels was supported by the DFG Research Center 
{\sc Matheon} in
Berlin.

 

\begin{thebibliography}{99}
\bibitem{CN}
{\sc B. D. Coleman, and  W. Noll}, {\em  
The thermodynamics of elastic materials with heat conduction and viscosity}, 
 Arch. Rational Mech. Anal. {\bf 13} (1963), pp. 167--178.

  
\bibitem{CGPS1} 
{\sc P.~Colli, G.~Gilardi, P.~Podio-Guidugli, and J.~Sprekels},
  {\em Existence and uniqueness of a global-in-time solution
  to a phase segregation problem of the Allen-Cahn type},
  Math. Models Methods Appl. Sci. {\bf 20} (2010), pp.~519--541.

\bibitem{CGPS2} 
{\sc P.~Colli, G.~Gilardi, P.~Podio-Guidugli, and J.~Sprekels}, 
{\em A temperature-dependent phase segregation problem of the
Allen-Cahn type}, Adv. Math. Sci. Appl. {\bf
20} (2010), pp.~219--234.

\bibitem{CGPS3} 
{\sc P.~Colli, G.~Gilardi, P.~Podio-Guidugli, and J.~Sprekels}, 
{\em Well-posedness and long-time behavior for a nonstandard viscous 
Cahn-Hilliard system}, WIAS Preprint No. 1602,
  Berlin 2011. Submitted. 

\bibitem{CGPS4}
{\sc P.~Colli, G.~Gilardi, P.~Podio-Guidugli, and J.~Sprekels},
{\em An asymptotic analysis for a nonstandard viscous 
Cahn-Hilliard system}. In preparation. 

\bibitem{FG}
{\sc E. Fried, and M. E. Gurtin},
{\em Continuum theory of thermally induced phase transitions based
on an order parameter}, Physica D {\bf 68} (1993), pp. 326--343.

\bibitem{Grie}
{\sc J. A. Griepentrog},
{\em Maximal regularity for nonsmooth parabolic problems
in Sobolev-Morrey spaces}, Adv. Differential Equations {\bf 12}
(2007), pp.~1031--1078. 
 
\bibitem{Gurtin}
{\sc M. Gurtin}, {\em Generalized Ginzburg-Landau and Cahn-Hilliard
equations based on a microforce balance}, Physica D {\bf 92},
pp. 178--192.

\bibitem{Hei} {\sc M. Heinkenschloss}, 
{\em The numerical solution of a control problem governed
by a phase field model}, Optim. Methods Softw. {\bf 7} (1997),
pp. 211--263.
 
\bibitem{HeiTr} {\sc M.~Heinkenschloss, and F. Tr\"oltzsch},
{\em Analysis of the Lagrange-SQP-Newton method for the
control of a phase field equation}, Control Cybernet.
{\bf 28} (1999), pp. 178-211.

\bibitem{HJ} {\sc K.-H. Hoffmann, and L. Jiang},
{\em Optimal control problem of a phase field model for solidification},
Numer. Funct. Anal. Optim. {\bf 40} (1992), pp. 11--27.
 
\bibitem{LSU} {\sc O. A. Lady\v{z}enskaya, V. A. Solonnikov,
and N. N. Ural'ceva}, {\em ``Linear and quasilinear equations of
parabolic type''}, Trans. Amer. Math. Soc. Vol. {\bf 23},
Providence, Rhode Island, 1968. 
 
\bibitem{LS07}
{\sc C.~Lefter, and J.~Sprekels},
{\em Control of a phase field system modeling non-isothermal
phase transitions}, Adv. Math. Sci. Appl. {\bf 17} (2007), pp.~181--194. 
    
\bibitem{Aubin}
{\sc J. L. Lions}, {\em ``Quelques m\'ethods de r\'esolution des
probl\`emes aux limites non lin\'eaires''}, Dunod Gauthier-Villars,
Paris, 1969.    
    
\bibitem{PG} {\sc P.~Podio-Guidugli}, {\em Models of phase segregation 
and diffusion of atomic species on a lattice}, Ric. Mat. 
{\bf 55} (2006), pp.~105--118.

\bibitem{Simon} {\sc J. Simon}, {\em Compact sets in the space}
$L^p(0,T;B)$, Ann. Mat. Pura. Appl. {\bf 146} (1987), pp.~65--96.

\bibitem{Tr} {\sc F. Tr\"oltzsch}, {\em ``Optimal Control of Partial
Differential Equations: Theory, Methods and Applications''}, Graduate Studies 
in Mathematics Vol. {\bf 112},
American Mathematical Society, Providence, Rhode Island, 2010. 
\end{thebibliography}
\end{document}